\documentclass[]{article}
\usepackage{lipsum}
\usepackage{amsfonts}
\usepackage{graphicx}
\usepackage{epstopdf}
\usepackage{lineno,hyperref}
\modulolinenumbers[5]
  \usepackage{paralist}
  \usepackage{hyperref}
  \usepackage{algorithm,algorithmic}
   
\ifpdf
  \DeclareGraphicsExtensions{.eps,.pdf,.png,.jpg}
\else
  \DeclareGraphicsExtensions{.eps}
\fi

\usepackage{graphicx} 
\usepackage{amsmath}
\usepackage{amssymb}
\usepackage{xcolor} 
\usepackage{hyperref}
\usepackage{tikz, pgfplots}
\usetikzlibrary{arrows,automata}

\usepackage{amssymb}
\usepackage {amsmath}
\usepackage{color,todonotes}
\usepackage{lineno,hyperref}
\modulolinenumbers[5]
  \usepackage{paralist}
  \usepackage{hyperref}  \usepackage{graphics} 
  \usepackage{epsfig}
\usepackage{graphicx}  
 \graphicspath{{./}{figures/}}
\usepackage{mathtools}
\DeclarePairedDelimiter{\ceil}{\lceil}{\rceil}

\DeclareMathOperator{\sgn}{sgn}

\newcommand{\I}{{\mathcal I}}
\newcommand{\J}{{\mathcal J}}

\newcommand{\RR}{{\mathbb R}}

\newcommand{\cQ}{{\mathcal Q}}

\allowdisplaybreaks

\newcommand{\R}{{\mathbb{R}}}

\newcommand{\cA}{{\mathcal A}}
\newcommand{\cB}{{\mathcal B}}

\newcommand{\cM}{{\mathcal M}}
\newcommand{\cN}{{\mathcal N}}

\newcommand{\cR}{{\mathcal R}}

\newcommand{\tr}{\operatorname{\text{tr}}}

\newcommand{\E}{{\mathbb E}}

\newtheorem{theorem}{Theorem}

\newtheorem{proposition}{Proposition}
\newtheorem{definition}{Definition}

\author{Simone Cacace \thanks{Dipartimento di Matematica e Fisica, Universit\`a Roma Tre, L.go S. Leonardo Murialdo, 1, 00146 Roma, Italy, 
                {\tt\small cacace@mat.uniroma3.it} }
 \and Roberto Ferretti
               \thanks{Dipartimento di Matematica e Fisica, Universit\`a Roma Tre, L.go S. Leonardo Murialdo, 1, 00146 Roma, Italy, 
                {\tt\small ferretti@mat.uniroma3.it} }
                \and Adriano Festa \thanks{Dipartimento di Ingegneria, Scienze dell'Informazione e Matematica, Universit\`a dell'Aquila, Via Vetoio, 67100 L'Aquila, Italy,
                {\tt\small adriano.festa@univaq.it} }
}

\title{\Large \bf Stochastic hybrid differential games\\and match race problems }

\begin{document}

\maketitle

\begin{abstract}
We discuss the general framework of a stochastic two-player, hybrid differential game, and we apply it to the modelling of a ``match race'' between two sailing boats, namely a competition in which the goal of both players is to proceed in the windward direction, while trying to slow down the other player. We provide a convergent approximation scheme for the computation of the value function of the game, and we validate the approach on some typical racing scenarios.
\end{abstract}

\noindent Stochastic hybrid systems, differential games, Hamilton--Jacobi equations 

\noindent
 93E20, 	49N70,	34K34,  65N06

\section{Introduction}

Hybrid processes are present in many economic and technological systems, whose dynamics can be modelled by a collection of controlled ordinary or stochastic differential equations: besides the standard actions performed on the current dynamics at a given time, the controller also has the option to switch to a different dynamics, in order to optimize some objective functional. 

Starting from the late 90s, several attempts have been made to provide a precise notion of hybrid systems. Among the different concepts proposed, we quote here \cite{branicky1998unified} and \cite{bensoussan1997hybrid} for respectively the deterministic and the stochastic case. The common feature of these models is to consider an extended state space for the dynamics, given by the product of both a continuous component and a discrete component, the latter indexing the active dynamics within a finite set.

In the optimal control of such class of systems, dynamic programming techniques have been widely investigated in the literature. The formulation of the problem in terms of a Bellman equation leads to a system of {\em quasi-variational inequalities}, which involve two different Bellman operators, related to respectively the continuous and the discrete control actions. A theoretical study of the problem in the framework of viscosity solutions can be found in \cite{bensoussan1997hybrid, dharmatti2005hybrid}. The numerical treatment via monotone schemes has also been studied in \cite{FerZid:2014}, proving that the classical Barles--Souganidis theory \cite{BarSou:1991} applies to the hybrid case, and providing a convergent technique to construct asymptotically optimal controls.

On the other hand, the case of differential games in the presence of hybrid dynamics seems much less explored in the literature. 
To our knowledge, the first study of a deterministic game under pure switching controls is given in   \cite{yong1990differential}.
Using the celebrated notion of non-anticipating strategies by Elliott and Kalton \cite{elliott1972existence}, one can prove the existence of a value for the game under a technical assumption, the so-called {\em no free loop} property. The unique value function satisfies a dynamic programming principle and can be characterized as the viscosity solution of an Isaacs system of quasi-variational inequalities.
We refer to this work also for an extensive review of the earlier literature on the subject. A more recent and general study, still in the deterministic case but much in the spirit of hybrid systems, is provided in \cite{shaiju2005differential}. Requiring also the classical \emph{Isaacs conditions}, the existence of a value is proved for games involving continuous, impulsive and discrete controls.  
Finally, concerning the stochastic case, we refer to some recent papers \cite{asri2018stochastic,ishii1991viscosity, hu2015switching, hu2010multi,hamadene2010switching, hamadene2013viscosity}. 

In this paper, we use the theory of stochastic hybrid differential games to model a route planning problem for two competing sailing boats, known as a {\em match race}. In this problem, the aim of the two competitors is to reach a mark at the end of a race leg \emph{before} the other, regardless of the time to reach the goal. This aspect makes the problem quite different from optimal navigation (discussed, e.g., in \cite{Spenkuch2011,dalang2015stochastic,ferrettifesta}) and motivates the use of game theory to model the interaction mechanics between the boats. This issue was also addressed in \cite{doi:10.1287/opre.1030.0078}, where the authors propose a technique to assess virtual competitions between yachts and to evaluate the pros and cons of various race scenarios. The physical interaction between the two boats (the \emph{wind shadow region} where one boat perturbs the wind) is modelled using a penalization/reward term.  Other works related are \cite{TAGLIAFERRI2014149,TAGLIAFERRI2017129}, where a risk model is included in the strategic decision process, and it is shown that, rather than finding the strategy that minimizes the time to complete the race leg, a strategy aimed at maximizing the probability of completing before the opponent offers better chances of victory. In the same works, the authors use a short term wind forecast methodology based on Artificial Neural Networks to model the instability of the wind.
The originality of the methodology that we propose in this paper lies on a game-theory-based formulation of the interaction between the two boats. This framework permits the observation of highly sophisticated strategic choices that are commonly used by tacticians in match race competitions, and a precise timing and quantification of them.

The paper is organized as follows. In Section \ref{Sect:hyb}, we introduce the mathematical framework for a stochastic hybrid game, reporting some results concerning its well-posedness in the viscosity sense, as well as the conditions for the existence of a value. In Section \ref{Sect:route}, we discuss in detail our game, and we analyze some relevant features of the corresponding value function. Section \ref{numerics} is devoted to the numerical solution of the Isaacs system of the game. We build a convergent algorithm based on a suitable monotone scheme, and we provide some hints on its actual implementation. Finally, in Section \ref{examples}, we perform some numerical tests, showing the effectiveness of the technique in different scenarios of application.

\section{Stochastic hybrid differential games} \label{Sect:hyb}


We describe the general structure of a zero-sum stochastic hybrid differential game, and we report the main results concerning the well-posedness of the problem in the sense of viscosity solutions.  The following presentation is a modified version of the one proposed in \cite{yong1990differential}. Fundamental contributions are also \cite{bensoussan1997hybrid, branicky1998unified} adapted to the stochastic case as in \cite{shaiju2005differential} or in the same spirit of  \cite{ferretti2014choosing}. We refer to these papers for further details and rigorous proofs. 

Given two compact sets $A\subset\R^{m_A}$  and $B\subset\R^{m_B}$ (for some integers $m_A$, $m_B$), we define the following standard sets of continuous controls for the two players, respectively 
$$
\mathcal{A}=\left\{a:(0,\infty) \to A \> | \> a \text{ measurable} \right\},$$
$$\mathcal{B}=\left\{b:(0,\infty) \to B \> | \> b \text{ measurable} \right\}.
$$
Moreover, in order to model the possibility for the two players to switch between different dynamics, we consider two finite sets of indices $\I=\{1,2,\ldots,N_\I\}$ and $\J=\{1,2,\ldots,N_\J\}$ (for some integers $N_\I$, $N_\J$), and we define the following sets of piecewise constant discrete controls, respectively
$$
\cQ = \left\{Q:(0,\infty) \to \I \> | \> Q(t)=\sum_{i\ge 0} q_i \chi_{[t_i,t_{i+1})}(t) \right\},$$
$$
\cR = \left\{R:(0,\infty) \to \J \> | \> R(t)=\sum_{i\ge 0} r_i \chi_{[t_i,t_{i+1})}(t) \right\},
$$
where $\{t_i\}$ is the sequence of (ordered) switching times, $\{q_i\}\subset\I$, $\{r_i\}\subset\J$ are the corresponding sequences of switching values for the two players, and $\chi_{[t_i,t_{i+1})}$ denotes the characteristic function of the interval $[t_i,t_{i+1})$.

We consider the dynamical system described by the following controlled stochastic differential equation (SDE):
\begin{equation} \label{eq_stato}
 \begin{cases}
 dX(t)=f(X(t),Q(t),a(t),R(t),b(t))dt+\sigma(X(t),Q(t),R(t))\,dW_t, \quad t>0\,,\\
 X(0)= x, \; Q(0^+)=q,\; R(0^+)=r,
 \end{cases}
\end{equation}
where (for some integers $d,k$) $x,X \in \R^d $, $q\in \I$, $r\in \J$, $a(\cdot)\in\mathcal{A}$, $Q(\cdot)\in\cQ$, $b(\cdot)\in\mathcal{B}$, $R(\cdot)\in\cR$, while $f :\R^d\times \I\times A \times \J \times B \to \R^d$ is the dynamics, $dW_t$ is the differential of a $k$-dimensional Brownian process, and 
$\sigma :\R^d\times \I \times \J \to \R^{d\times k}$ is the corresponding covariance matrix.

To properly define solutions of the stochastic differential equation \eqref{eq_stato}, we need a standard regularity assumption:
\begin{itemize}
 \item[{\bf H1 -}]$f$ and $\sigma$ are globally bounded and uniformly Lipschitz continuous with respect to $x$.
\end{itemize}
Then, the following integral representation formula holds:
$$ X(t)= x+\int_0^t f(X(s), Q(s), a(s),R(s), b(s))ds+\int_0^t\sigma(X(s), Q(s),\bar R(s))\,dW_s.$$
The stochastic trajectory starts from $(x,q,r)$ in the extended state space $\R^d\times \I \times \J$. At each time $t>0$ the first player can act on the current dynamics through the control $a(\cdot)\in\mathcal{A}$, or switch to another dynamics using the discrete control $Q(\cdot)\in\cQ$. Similarly, the second player employs the controls $b(\cdot)\in\mathcal{B}$ and $R(\cdot)\in\cR$. 
This setting is suitable for our application to a match race competition, but we remark that the most general framework of hybrid control systems (see, e.g., \cite{bensoussan1997hybrid}), allows one to deal with problems including also autonomous transitions and jumps in the state $X$.

Now, we define the game between the two players. To this end, let us introduce a more compact notation for the controllers, by setting respectively $\alpha(t):=(Q(t),a(t))\in \cQ\times \cA$ and $\beta(t):=(R(t),b(t))\in \cR\times \cB$. Moreover, we consider  the following cost functional:
\begin{multline}\label{J}
 J(x,q,r;\alpha,\beta)  := \E \left(\int_0^{+\infty} e^{-\lambda s} \ell(X(s),Q(s),a(s),R(s),b(s))ds \right.\\ \left.+ \sum_{i\ge 0} e^{-\lambda  t_i} \left[C_A\left(Q( t_i^-),Q( t_i^+)\right)+C_B\left(R( t_i^-),R( t_i^+)\right)\right]\right)\,.
\end{multline}
Here, the symbol $\E$ denotes expectation with respect to the Wiener measure, while the first integral term defines a standard infinite horizon functional, with discount factor $\lambda>0$  and a running cost $\ell:\R^d\times\cQ\times \cA\times \cR\times \cB\to\R$. We assume that:
\begin{itemize}
 \item[{\bf H2 -}] $\ell$ is non-negative, bounded and uniformly Lipschitz continuous with respect to $x$.
\end{itemize}
On the other hand, the second term in \eqref{J} accounts for the discounted costs $C_A:\I\times\I\to\R$ and $C_B:\J\times\J\to\R$ associated to the switches of the two players ($A$ and $B$ respectively) at times $\{t_i\}$. Here, player $A$ wants to maximize $J$ using the control $\alpha$, thus paying a {\em negative} cost $C_A$ for each switch. Similarly, player $B$ wants to minimize $J$ using the control $\beta$ and paying a {\em positive} cost $C_B$ for each switch. Note that, to simplify notation, we regrouped the switching times of both players in a single sequence $\{t_i\}$. This means that, if {\em only} one player performs a switch at time $t_i$, the corresponding cost of the other player should be zero. We summarize all these properties by requiring the following assumptions:
\begin{itemize}
 \item[{\bf H3 -}] $C_A$ and $C_B$ are bounded and satisfy  
  $$
  C_A(q,q)=0\quad \mbox{for every }q\in\I\,,
  \qquad
  C_B(r,r)=0\quad \mbox{for every }r\in\J\,.
  $$
  Moreover, there exists $C_0>0$ such that
  $$
  \max_{q_1\neq q_2}C_A(q_1,q_2)\leq -C_0, \qquad  \min_{r_1\neq r_2}C_B(r_1,r_2)\geq C_0.
  $$
\end{itemize}

We proceed by defining the value functions of the game. To this end, we employ the classical notion of \emph{non-anticipating strategies} \cite{elliott1972existence,yong1990differential}, which allows to rigorously prove a dynamic programming principle. 
\begin{definition}
    A \emph{non-anticipating strategy} for player $A$ (resp. $B$) is a map $\phi:\cR\times\mathcal{B}\rightarrow \cQ\times\mathcal{A}$ (resp. $\psi:\cQ\times\mathcal{A}\rightarrow \cR\times\mathcal{B}$) such that,  for any $t>0$,  
    $$\beta(s)=\tilde{\beta}(s)\mbox{\, for all \,}s\leq t\mbox{\, implies \,} \phi[\beta](s)=\phi[\tilde{\beta}](s)\mbox{\, for all \,}s\leq t\,.$$ 
    $$\mbox{(resp. \,}\alpha(s)=\tilde{\alpha}(s)\mbox{\, for all \,}s\leq t\mbox{\, implies \,}\psi[\alpha](s)=\psi[\tilde{\alpha}](s)\mbox{\, for all \,}s\leq t\,.)$$ 
\end{definition}

We denote the set of non-anticipating strategies by $\Phi$ for player $A$, and by $\Psi$ for player $B$.
Then, for every $(x,q,r)\in\R^d\times \I \times \J$, we define the \emph{lower value} function $\underline{v}$ of the game as
\begin{eqnarray} \label{l_value}
 \underline{v}(x,q,r):= \inf_{\psi\in \Psi} \sup_{\alpha\in\cQ\times\cA}J(x,q,r;\alpha,\psi[\alpha]),
\end{eqnarray}
and the \emph{upper value} $\overline{v}$ as 
\begin{eqnarray} \label{u_value}
 \overline{v}(x,q,r):= \sup_{\phi\in \Phi} \inf_{\beta\in\cR\times\cB}J(x,q,r;\phi[\beta],\beta).
\end{eqnarray}
Moreover, if $ \underline{v}\equiv \overline{v}$, we say that \emph{the game has a value}, and we denote it by $v$.

In the next Proposition, we state the dynamic programming principle satisfied by both the value functions.
\begin{proposition}\label{PropDPP}Under the assumptions H1-H3, for all $(x,q,r) \in \R^d\times\I\times\J$ and $\tau > 0$, the following equation holds true
\begin{multline}
\label{ppdgam}
\underline{v}(x,q,r)=\inf_{\psi\in \Psi} \sup_{\alpha\in\cQ\times\cA}\Big\{ \E\Big(\int_0^\tau\ell(X(s),\alpha(s),\psi[\alpha](s))ds \\
+ \sum_{t_i<\tau} e^{-\lambda  t_i} \left[ C_A\left(Q( t_i^-),Q( t_i^+)\right)+ C_B\left(R( t_i^-),R( t_i^+)\right)\right]\\
+\underline{v}(X(\tau), Q(\tau), R(\tau))e^{-\lambda \tau}\Big)\Big\},
\end{multline}
where $Q$ and $R$ are the switching controls contained respectively in the strategy $\alpha$ and $\psi[\alpha]$.
A similar equation holds for the upper value function $\overline{v}$, by swapping the role between $\inf$ and $\sup$ in \eqref{ppdgam}.
\end{proposition}
Now, for a generic function $\varphi:\R^d\times \I\times \J\to \R$, we define the two following {\em switching operators} 
\begin{eqnarray*}
 \cN [\varphi](x,q,r) := \max_{\hat q\neq q} \{\varphi(x,\hat q,r)+ C_A(q,\hat q)\},\\
 \cM [\varphi](x,q,r) := \min_{\hat r\neq r} \{\varphi(x,q,\hat r)+ C_B(r,\hat r)\},
\end{eqnarray*}
which provide some natural bounds on the value functions, as stated in the next Proposition.
\begin{proposition}\label{lem}
For every $(x,q,r)\in\R^d\times\I\times\J$ the lower value function $\underline v$ satisfies
$$ \cN[\underline v] (x,q,r)   \leq \underline v(x,q,r) \leq  \cM[\underline v](x,q,r)\,.$$ The same estimates hold for the upper value function $\overline v$.
\end{proposition}

Proposition \ref{PropDPP} and Proposition \ref{lem} allow to derive the Hamilton--Jacobi--Isaacs equations of the game. More precisely, for $x,p\in\RR^d$, $q\in\I$ and $r\in\J$, we introduce the Hamiltonians 
\begin{eqnarray}\label{Ham}
 H^-(x,q,r,p) := \min_{a\in A}\max_{b\in B}\{ - f(x,q,a,r,b)\cdot p - \ell(x,q,a,r,b)\}, \\ H^+(x,q,r,p) :=  \max_{b\in B} \min_{a\in A}\{ - f(x,q,a,r,b)\cdot p - \ell(x,q,a,r,b)\},
\end{eqnarray}
and the second order differential operators
\begin{eqnarray}\label{Diffop}
 F^-[\varphi](x,q,r)=\lambda \varphi(x,q,r)+ H^-(x,q,r,D \varphi)-\frac{1}{2}\tr\left( \sigma\sigma^T D^2 \varphi(x,q,r)\right), \\ F^+[\varphi](x,q,r)=\lambda \varphi(x,q,r)+ H^+(x,q,r,D \varphi)-\frac{1}{2}\tr\left( \sigma\sigma^T D^2 \varphi(x,q,r)\right),
\end{eqnarray}
where $D$ and $D^2$ denote respectively the gradient and the hessian with respect to $x$,  $\sigma^T$ is the transpose of $\sigma$, and $\tr(\cdot)$ stands for the matrix trace.

Then, it follows that the value functions $\underline v$ and $\overline v$ satisfy, for every $(x,q,r)\in\R^d\times\I\times\J$, respectively
\begin{equation}\label{hjb1}
\max\left\{\underline v-\cM [\underline v], \min\left\{\underline v-\cN [\underline v],F^-[\underline v] \right\}\right\}= 0\,,
\end{equation}
and 
\begin{equation}\label{hjb2}
\max\left\{\overline v-\cM [\overline v],  \min\left\{\overline v-\cN [\overline v],F^+[\overline v]\right\}\right\}= 0\,,
\end{equation}
namely two systems of $N_\I N_\J$ quasi-variational inequalities. In each system, we can identify three separate operators, which provide respectively the best possible switching for the two players, and the best possible continuous controls. The arguments attaining the respective extrema in such equations represent the overall optimal control strategies. The derivation of \eqref{hjb1} and \eqref{hjb2}, which is elementary under differentiability assumptions, can be rigorously justified in a more general setting by an adaptation of the viscosity theory \cite{dharmatti2005hybrid} to the case under consideration. 

To conclude this section, we briefly discuss the key steps for proving the existence of a value for the game, namely that $\underline v\equiv\overline v$. First, a uniqueness result for the viscosity solutions of both \eqref{hjb1} and \eqref{hjb2} is needed. 
In this direction, the following additional assumption, the so-called {\em no free loop} property, appears in several papers on hybrid games, see e.g., \cite{hamadene2013viscosity,yong1990differential,ishii1991viscosity}:
\begin{itemize}
 \item[{\bf H4 -}]Let $\{(q_i,r_i)\}_{i=1,...,N+1}$ be a finite sequence of indices such that $(q_i,r_i)\not =(q_{i+1},r_{i+1})$ for every $i=1,...,N-1$ and  $(q_1,r_1)=(q_{N+1},r_{N+1})$. Then $$  \sum_{i=1}^{N} \left\{ C_A(q_i,q_{i+1})+ C_B(r_i,r_{i+1})\right\}\neq 0.$$
\end{itemize}
Although technical, this assumption seems unavoidable in order to obtain a comparison principle between a viscosity sub-solution $u$ and a viscosity super-solution $w$ of \eqref{hjb1} (the same reasoning applies to \eqref{hjb2}). The idea is that, using assumption H4, one can find, for every $x\in\R^d$, a common state $(q^*,r^*)\in\I\times\J$ in which both inequalities $F^-[u](x,q^*,r^*)\le 0$ and $F^-[w](x,q^*,r^*)\ge 0$ hold. Then, one can proceed with the usual comparison of the Hamiltonians in the viscosity theory and conclude that $u\le w$. 
This result implies that $\underline v$ is the unique viscosity solution of \eqref{hjb1} and $\overline v$ is the unique viscosity solution of \eqref{hjb2}.

Finally, the existence of a value for the game is guaranteed by providing assumptions that let the Isaacs systems \eqref{hjb1} and \eqref{hjb2} coincide, as for the following classical {\em Isaacs conditions}:
\begin{itemize}
\item[{\bf H5 -}] 
$H^-(x,q,r,p)=H^+(x,q,r,p)$ for every $(x,q,r,p)\in\R^d\times\I\times\J\times\R^d$.
\end{itemize}
Summarizing, we have the following result.
\begin{theorem}\label{gen-isaacs}
Under assumption H1-H5, the value function $v:=\underline{v}\equiv\overline{v}$ is the unique viscosity solution of both \eqref{hjb1} and \eqref{hjb2}.
\end{theorem}

\section{The match race problem}\label{Sect:route}

We apply the theoretical framework of hybrid differential games, discussed in the previous section, to a real-world application. A {\em match race} is a competition between two sailing boats, in which the goal of both players is to reach, as first, the end of the course, regardless of their relative distance. Each player can take advantage of the wind fluctuations to proceed upwind towards the finish line, by adjusting the relative angle between the sail and the wind, and also changing the tack side. In addition, the players can make use of their respective influence, caused by the sail turbulence. This is usually an advantage for the leading boat, which can exploit this influence to control the other player. Note that a reasonable description of this problem requires, at least, a state space of dimension $d=5$: two pairs of coordinates to track the positions of the boats in a plane, and one coordinate for the wind angle. Here, we consider a simplified game, namely we neglect the windward mark, and we just focus on the strategies of the two players in the space of relative positions. This reduces the problem to dimension $d=3$, but it is still a realistic racing criterion when the two players are far from the windward mark. Moreover, it can be numerically solved in a reasonable time also on a laptop computer. The analysis and the parallel implementation of the full game is under investigation, and will be addressed in a forthcoming paper.

In the next subsections, we first introduce the hybrid dynamics for the boats, following the model presented in \cite{ferrettifesta}. Then, we define the hybrid game, by suitably setting all the parameters appearing in the cost functional \eqref{J}. Finally, we present a more detailed analysis of the value function of the game, in case the two players are far enough from each other.

\subsection{Dynamics modelling}\label{dynmod}
\begin{figure}[!h]
\centering
\begin{tabular}{ccc}
\includegraphics[width=.31\textwidth]{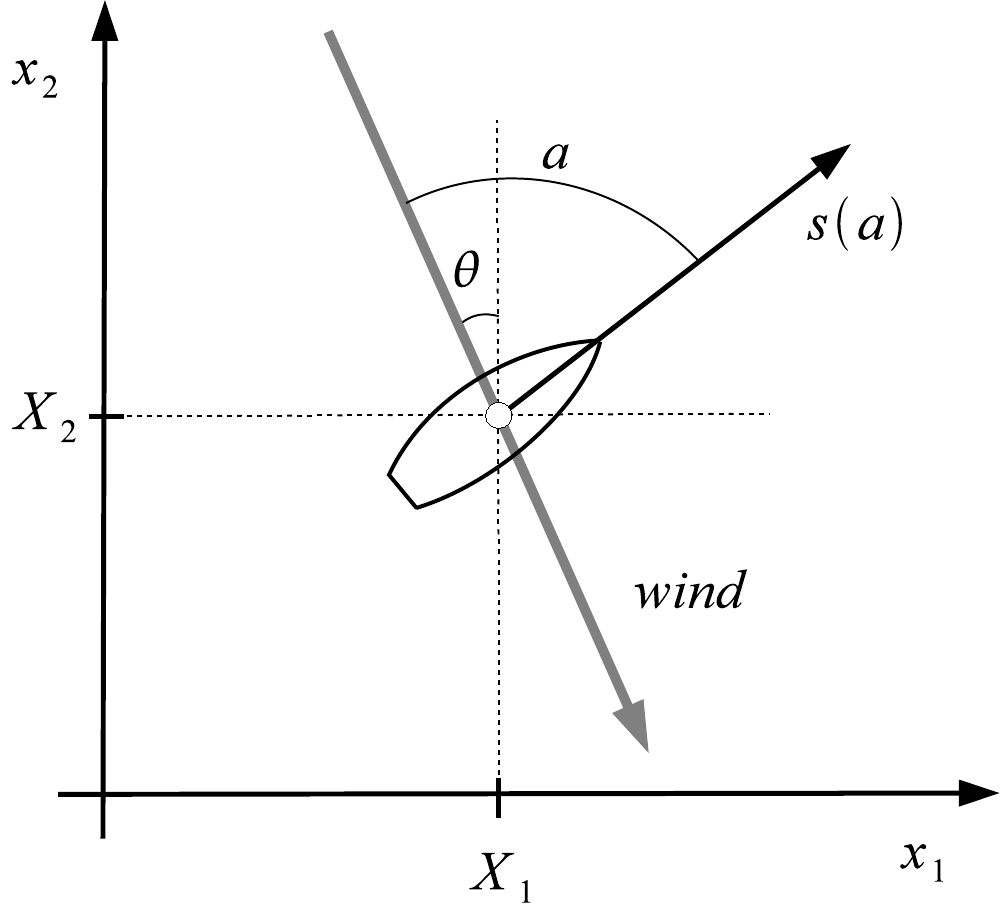} &
\includegraphics[width=.31\textwidth]{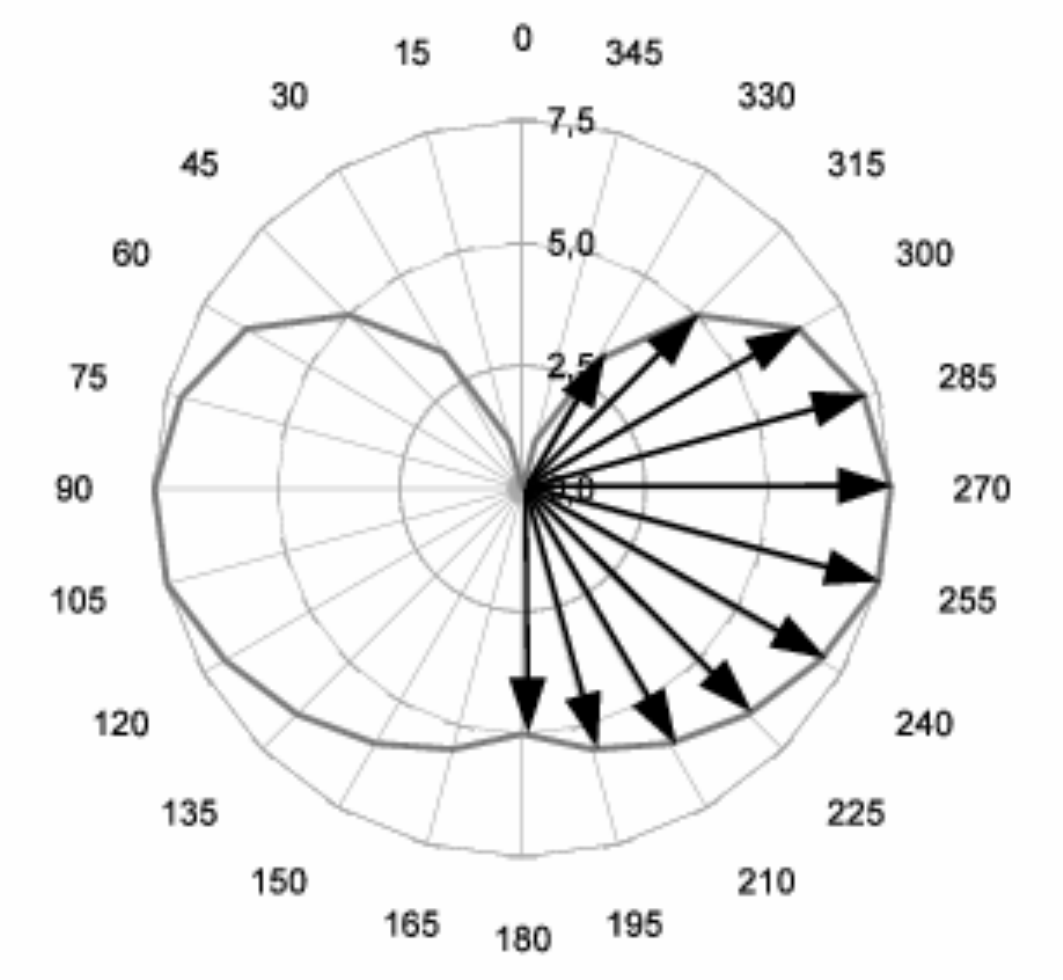} &
\includegraphics[width=.31\textwidth]{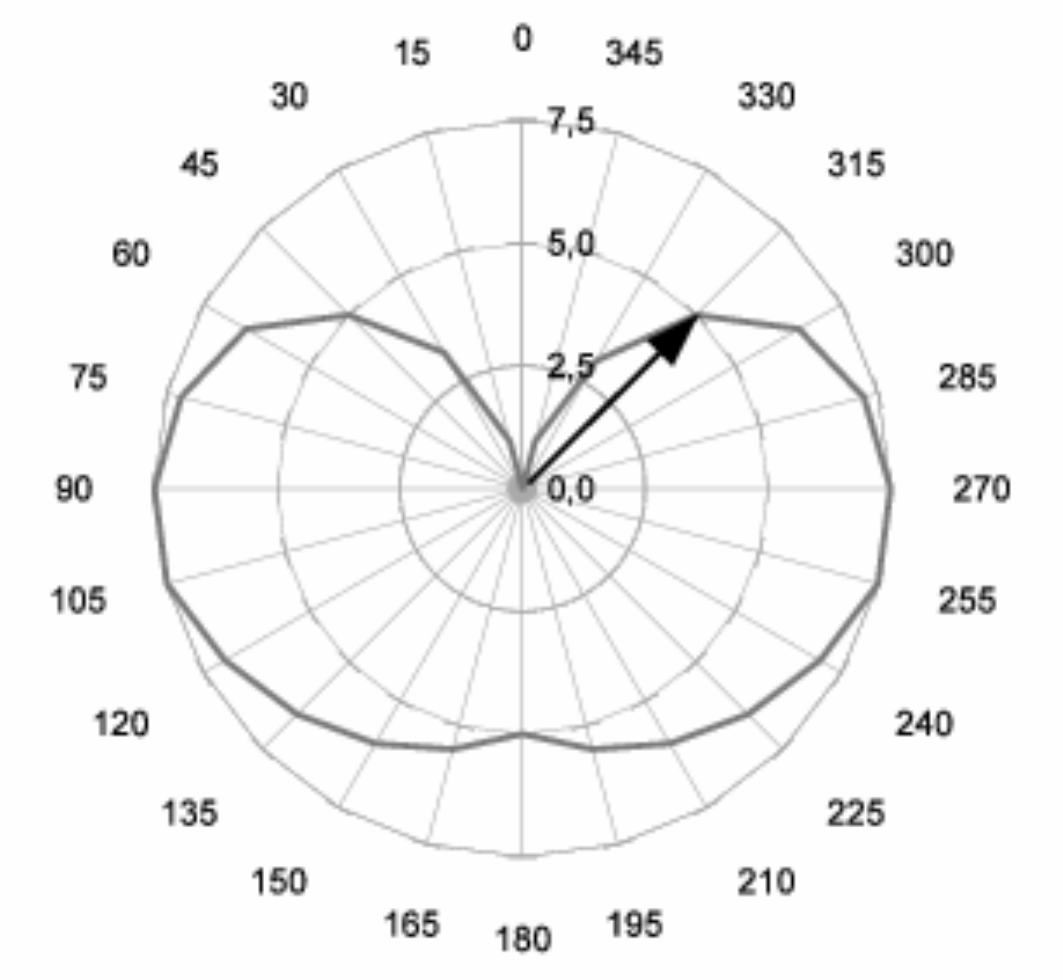} \\
(a)&(b)&(c)
\end{tabular}
\caption{Model of the boat speed. Geometric setting (a), one of the two dynamics ($q=1$) of a boat, superposed on the polar plot of the speed (b) and simplified dynamics based on the angle of largest windward component of the speed (c).} \label{fig:dyn}
\end{figure}

We consider the motion on a plane of a single boat, subject to a wind of constant speed and variable direction. We set the dimension of the state space to $d=3$, in which the first two components $x_1$ and $x_2$ represent the position of the boat, while the third component $x_3$ gives the angle $\theta\in[-\pi,\pi]$ of the wind with respect to the vertical axis. In particular, $\theta$ is negative in the second and third quadrant, and positive in the first and fourth quadrant, see Figure \ref{fig:dyn}a. 

Moreover, we assume that the wind has a purely Brownian nature, i.e., it evolves according to the one-dimensional SDE:
\begin{equation}\label{eq:simpl_wind}
dX_3(t)=d\Theta(t) = \sigma dW(t),
\end{equation}
where $dW$ denotes the differential of a Brownian process, and $\sigma>0$ is the corresponding standard deviation.

On the other hand, the motion of the boat results from both the wind direction and the boat characteristics. Following the notation presented in Section \ref{Sect:hyb}, we introduce the set of admissible controls $A=[0,\pi]$ as the unsigned angles between the boat direction and the wind, so that the continuous control is given by a function $a:[0,+\infty)\to A$. Then, since the wind speed is constant, the boat speed will depend only on the angle $a$, by means of a function $s:A \rightarrow \R_+$, the so-called \emph{polar plot} of the boat. Figure \ref{fig:dyn}a summarizes this geometric setting, while in Figure \ref{fig:dyn}b we show a typical form of the polar plot, with the whole set of speeds associated with the port tack. 
Note that, for $a=0$, the trajectory points directly in the upwind direction, whereas, for $a=\pi$, the trajectory has the same direction of the wind field.

When sailing to windward, it is customary to keep constantly the boat at its most efficient angle with the wind, that is, at the angle $a^*\approx \frac{\pi}{4}$, corresponding the largest windward component of the speed. In this case, the dynamics can be simplified by freezing the control at the value $a^*$ (see Figure \ref{fig:dyn}c), and acting on the system only by changing tack. In what follows, we will use this simplified approach.

Finally, we introduce the discrete control, namely a piecewise constant function $Q:[0,+\infty)\to \I$, taking values in the discrete set $\I=\{1,2\}$. The two possible discrete states correspond to the tack sides, where the port tack is identified by $q=1$ and the starboard tack by $q=2$. Hence, the dynamics of the boat is given by
\begin{equation}\label{eq:simpl_dyn}
\begin{cases}
\dot X_1(t) = s(a^*)\sin\left(\Theta(t) + (-1)^{Q(t)} a^*\right) \\
\dot X_2(t) = s(a^*)\cos\left(\Theta(t) + (-1)^{Q(t)} a^*\right).
 \end{cases}
\end{equation}
\subsection{Game modelling}
\begin{figure}[!h]
\begin{center}
\begin{tabular}{c}
\includegraphics[width=.5\textwidth]{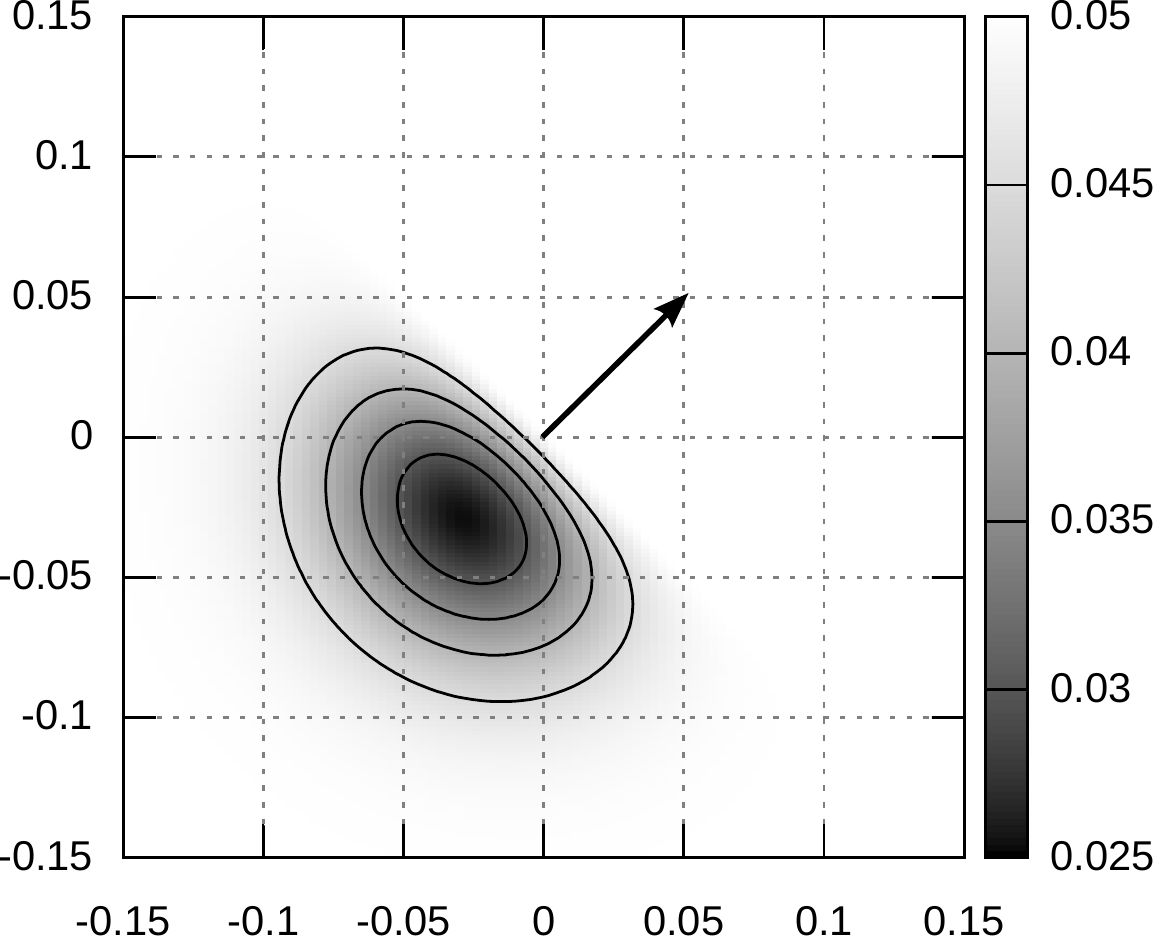}
\end{tabular}
\end{center}
\caption{Level sets of the speed profile.}\label{speed_profile}
\end{figure}

We define the game in \emph{reduced coordinates}, i.e., we consider as space variable the relative position of the two players. We denote by $x^A=(x_1^A,x_2^A)\in\R^2$, $x^B=(x_1^B,x_2^B)\in\R^2$ and $\theta\in[-\pi,\pi]$, respectively the coordinates of the two players and the wind angle, while the reduced coordinates are given by $x=x^A-x^B\in\R^2$. Then, for 
$q,r\in\mathcal{I}=\mathcal{J}=\{1,2\}$ and discrete controls $Q,R:[0,+\infty]\to\I$ such that $Q(0)=q$ and $R(0)=r$, we define 
the controlled dynamics of the game according to \eqref{eq:simpl_wind} and \eqref{eq:simpl_dyn}:
\begin{equation}\label{game-dyn}
\left\{\begin{array}{l}
        dX^A(t)=f^A(X(t),\Theta(t),Q(t))dt\\
        dX^B(t)=f^B(X(t),\Theta(t),R(t))dt\\
        d\Theta(t)=\sigma dW(t)
       \end{array}
\right.
\qquad
\left\{\begin{array}{l}
        X^A(0)=x^A\\
        X^B(0)=x^B\\
        \Theta(0)=\theta
       \end{array}
\right.
\end{equation}
where 
$$
f^A(x,\theta,q)=s^A(x,\theta)\left(\sin(\theta+(-1)^{q} a^*),\cos(\theta+(-1)^{q} a^*)\right),
$$
$$
f^B(x,\theta,r)=s^B(-x,\theta)\left(\sin(\theta+(-1)^{r} b^*),\cos(\theta+(-1)^{r} b^*)\right),
$$
with $a^*=b^*=\frac{\pi}{4}$. The speed functions $s^A$ and $s^B$ contain the information about the interaction between the two players.
As an example, we can take
\begin{equation}\label{speedex}
s^P(x,\theta)=\bar s^P\left(1+\min\{s_0^P(x\cdot (\sin(\theta),\cos(\theta))e^{-s_1^P|x|^2},0\}\right) \quad (P=A,B),
\end{equation}
for given positive constants $\bar s^P, s_0^P, s_1^P$, which would model the wind shadow region, i.e., a situation in which the player $P$ has its maximum speed $\bar s^P$ when the two players are far from one another, but it is slowed down when its position is close and behind or on the downwind side of the other (note the dependency of $f^A$ on $x$ and of $f^B$ on $-x$, which reflects the speed profile with respect to the origin, according to the leading player). Figure \ref{speed_profile} shows the level sets of the speed function $s^P$ in \eqref{speedex}, corresponding to $\bar s^P=0.05$, $s_0^P=20$ and $s_1^P=300$, for $\theta=\frac{\pi}{4}$. 

We remark that we can also incorporate in the maximum speeds $\bar s^A$ and $\bar s^B$ an additional dependency on the states $(x,q,r)$. This can be useful to model the rules about the right of way in match race competitions. For instance, with a little abuse of notation, we can choose 
$$\bar s^A(x,\theta,q)=\bar s^A\left(1-\nu_1  e^{-\frac{x^2}{\nu_2}}\right)^{q-1}$$
(and similarly $\bar s^B(x,\theta,r)$) to introduce a penalization of the speed of the boat on the \emph{port-tack} ($q,r=2$),  which activates only when the distance between the boats is small enough, tuned by the choice of the parameters $\nu_1, \nu_2$. This penalization mimics the interaction of a boat meeting a competitor on opposite tacks (see \cite[Rule 10]{isaf}).

In order to define the infinite horizon cost functional $J$ in \eqref{J}, we choose two constants $C^A, C^B>0$, and we set the switching costs as follows
$$
C_A(q_1,q_2)=\left\{
\begin{array}{ll}
    -C^A & \mbox{if } q_1\neq q_2 \\
    0 & \mbox{otherwise,} 
\end{array}\right.
\qquad C_B(r_1,r_2)=\left\{
\begin{array}{ll}
    C^B & \mbox{if } r_1\neq r_2 \\
    0 & \mbox{otherwise.} 
\end{array}\right.
$$
Finally, we choose the running cost
$$
\ell(x,\theta,q,r)=f^A_2(x,\theta,q)-f^B_2(x,\theta,r),
$$
so that the cost functional integrates the vertical component of the relative speed of the two boats. 
This results in a game in which  each player wants to overcome the opponent along the vertical component with the least number of switches. Due to the lack of continuous controls, the corresponding systems of quasi-variational inequalities \eqref{hjb1}-\eqref{hjb2} coincide, and take the form:
\begin{equation}\label{mr-qvi}
\begin{split}
\min \Big\{& v(x,\theta,q,r)-v(x,\theta,\hat q,r)+C^A\,,\, \max\Big\{v(x,\theta,q,r)-v(x,\theta,q,\hat r)-C^B\,,\\& \lambda v(x,\theta,q,r)- f(x,\theta,q,r)\cdot D v - \ell(x,\theta,q,r) -\frac{\sigma^2}{2}\frac{\partial^2 v}{\partial \theta^2} (x,\theta,q,r) \Big\}\Big\}= 0\,,
\end{split}
\end{equation}
where, for every $q,r\in\I=\J=\{1,2\}$, we set $\hat q=3-q$ and $\hat r=3-r$, and we denoted by $f$ the deterministic part of the coupled dynamics in $\R^d$, namely
$$
f(x,\theta,q,r)=\left(f^A_1(x,\theta,q)-f^B_1(x,\theta,r),f^A_2(x,\theta,q)-f^B_2(x,\theta,r), 0\right).
$$
We can observe that, in the present setting, the technical assumption H4 in Theorem \ref{gen-isaacs} is satisfied if $C_A\neq C_B$. Otherwise, uniqueness of a solution is not ensured. 

\subsection{Decoupling of the game with ``far'' players}
As discussed in the previous section, the coupling in the dynamics of the two players, and hence the essence of the game, is entirely embedded in the speed function $s^P$. 
A key observation is that if the two players are far enough from each other, i.e., if $|x|\gg1$, then $s^P\approx \bar s^P$.  Consequently, a fair approximation of the {\em far dynamics} of each player depends only on the wind direction and on the switching strategy. 
In this setting, we can provide a more explicit analysis of the game, and also obtain suitable boundary conditions for the approximation of the problem in a bounded domain, as it will be discussed later. We remark that this analysis is much in the same spirit of the one carried out in \cite[Chapter 5]{vinck}, and in some sense brings it to its final conclusions, in the case in which the player is far from the target.

Assuming that $|x|=|x^A-x^B|\gg1$, and using the definition of $s^P$ and $\ell$, we can split the cost functional $J$ in \eqref{J} as the difference 
$$
J(x,\theta,q,r;Q,R)=J^{A}(\theta,q;Q)-J^{B}(\theta,r;R),
$$
where
$$
J^A(\theta,q;Q)=\mathbb{E}\left(\int_0^\infty \bar s^A \cos\left(\Theta(t)+\frac{\pi}{4}(-1)^{Q(t)}\right)e^{-\lambda t}dt-C^A\sum_{i\ge 0}e^{-\lambda t_i^A}\right),
$$
$$
J^B(\theta,r;R)=\mathbb{E}\left(\int_0^\infty \bar s^B\cos\left(\Theta(t)+\frac{\pi}{4}(-1)^{R(t)}\right)e^{-\lambda t}dt-C^B\sum_{i\ge 0}e^{-\lambda t_i^B}\right).
$$
As a consequence, we  get
\begin{align*}
v(x,\theta,q,r)&=\inf_{R(\cdot)}\sup_{Q(\cdot)}J(x,\theta,q,r;Q,R)\\
                    &=\sup_{Q(\cdot)}J^{A}(\theta,q;Q)+\inf_{R(\cdot)}\left\{-J^{B}(\theta,r;R)\right\}\\
                    &=\sup_{Q(\cdot)}J^{A}(\theta,q;Q)-\sup_{R(\cdot)}J^{B}(\theta,r;R)\\
                    &=v^{A}(q,\theta)-v^{B}(r,\theta),
\end{align*}
where, for $P=A,B$, and $p=q,r\in\I=\J=\{1,2\}$, we denote by $v^P(p,\theta)$ the value function corresponding to the optimal control problem, 
for the single player $P$, of maximizing $J^P$ subject to the dynamics $f^P$. We remark that, due to the special structure of $f^P$ and of the running cost in $J^P$, the value function $v^P$ 
depends only on $\theta$ and on the discrete state $p$. Moreover, it satisfies the following system of quasi-variational inequalities: for $p\in\I=\{1,2\}$ and $\hat p=3-p$,
\begin{equation}\label{singleplayerQVI-theta}
\begin{split}
\min \Big( & v^P(p,\theta)-v^P(\hat p,\theta)+C^P ,\\
           &\lambda v^P(p,\theta)-\bar s^P\cos\left(\theta+\frac{\pi}{4}(-1)^p\right)-\frac{\sigma^2}{2} \frac{\partial^2 v^P}{\partial \theta^2} (p,\theta)\Big)=0.
\end{split}
\end{equation}
For general switching costs $C^A$, $C^B$ and speeds $\bar s^A$, $\bar s^B$, we can solve \eqref{singleplayerQVI-theta} numerically, as shown in the next section. Neverthless, in the symmetric case ($C^A=C^B=:\bar c$ and $\bar s^A=\bar s^B=:\bar s$), we have $v^A(p,\theta)=v^B(p,\theta)=:\bar v(p,\theta)$, and we can extract further information by straightforward computations. 
Indeed, choosing alternately $p=1$ and $p=2$ in \eqref{singleplayerQVI-theta}, 
for every $\theta\in[-\pi,\pi]$ we get 
$$
\bar v(1,\theta)-\bar v(2,\theta)+\bar c\ge 0,\qquad \bar v(2,\theta)-\bar v(1,\theta)+\bar c\ge 0,
$$
and, at the points $\theta$ such that both inequalities are strict, we also have the equations
$$
\lambda \bar v(1,\theta)-\bar s\cos\left(\theta-\frac{\pi}{4}\right)-\frac{\sigma^2}{2} \frac{\partial^2 \bar v}{\partial \theta^2} (1,\theta)=0\,,
$$
$$
\lambda \bar v(2,\theta)-\bar s\cos\left(\theta+\frac{\pi}{4}\right)-\frac{\sigma^2}{2} \frac{\partial^2 \bar v}{\partial \theta^2} (2,\theta)=0\,.
$$
Defining the difference $\tilde v(\theta)=\bar v(2,\theta)-\bar v(1,\theta)$, by linearity we readily obtain
\begin{equation}\label{double-obstacle}
\max\left\{\tilde v(\theta)-\bar c\,,\,\min\left\{ \tilde v(\theta)+\bar c\,,\,
 \lambda \tilde v(\theta)-\sqrt{2}\bar s\sin(\theta)-\frac12 \sigma^2 \frac{\partial^2}{\partial\theta^2}\tilde v(\theta)\right\}\right\}=0,
\end{equation}
where we used the subtraction formula for the cosine function. This is a classical double obstacle problem, whose solution can be characterized as follows. First of all, the solution to the second order differential equation is given by
$$
v^*(\theta)=C_1e^{-\omega^*\theta}+C_2e^{\omega^*\theta}+\Omega^*\sin(\theta),\quad\omega^*=\displaystyle\frac{\sqrt{2\lambda}}{\sigma}\,,\quad  \Omega^*=\displaystyle\frac{2\sqrt2 \bar s}{2\lambda+\sigma^2}\,,$$ where $C_1, C_2$ are constants to be determined. 
By symmetry we require $v^*(0)=0$, whereas, imposing $C^1$ regularity for the contact point $\theta^*$ with the obstacle (this is a classical result, see \cite{MR0318650}), we get $v^*(\theta^*)=\bar c$ and $\frac{\partial v^*}{\partial \theta}(\theta^*)=0$. 
This easily implies the following nonlinear equation in $\theta$:
\begin{equation}\label{opt-comm}
\Omega^*\sin(\theta)-\frac{\Omega^*}{\omega^*}\tanh(\omega^*\theta)\cos(\theta)=\bar c\,,
\end{equation}
which admits a unique solution $\theta^*\in[0,\frac{\pi}{2}]$, since the left hand side is strictly increasing for $\theta\in[0,\frac{\pi}{2}]$ (we recall that, when sailing to windward, this is the interesting case). 

Hence, we obtain 
\begin{equation}\label{opt-val}
\tilde v(\theta)=\left\{\begin{array}{ll}
                    -\bar c & \theta<-\theta^*\\
                     v^*(\theta) & |\theta|\le \theta^*\\
                    \bar c & \theta>\theta^*
                   \end{array}
\right.
\end{equation}
and coming back to the relationship
$$
v(\theta,q,r)=v^{A}(\theta,q)-v^{B}(\theta,r)\,,
$$
we conclude that, for $|x|\gg 1$ 
$$
v(x,\theta,1,1)=\bar v(1,\theta)-\bar v(1,\theta)=0\,,\qquad v(x,\theta,2,2)=v(2,\theta)-\bar v(2,\theta)=0\,,
$$
$$
v(x,\theta,2,1)=\bar v(2,\theta)-\bar v(1,\theta)=\tilde v(\theta)\,,\qquad v(x,\theta,1,2)=\bar v(1,\theta)-\bar v(2,\theta)=-\tilde v(\theta)\,. 
$$
As a final remark, we point out again that, if the game is not symmetric ($C^A\not = C^B$ or $\bar s^A\not=\bar s^B$), then no such explicit computation is possible, and the single-player solution, as well as the boundary conditions for the two-player game, should be computed numerically.

\begin{figure}[!t]
\begin{center}
\begin{tabular}{c}
\includegraphics[width=.5\textwidth]{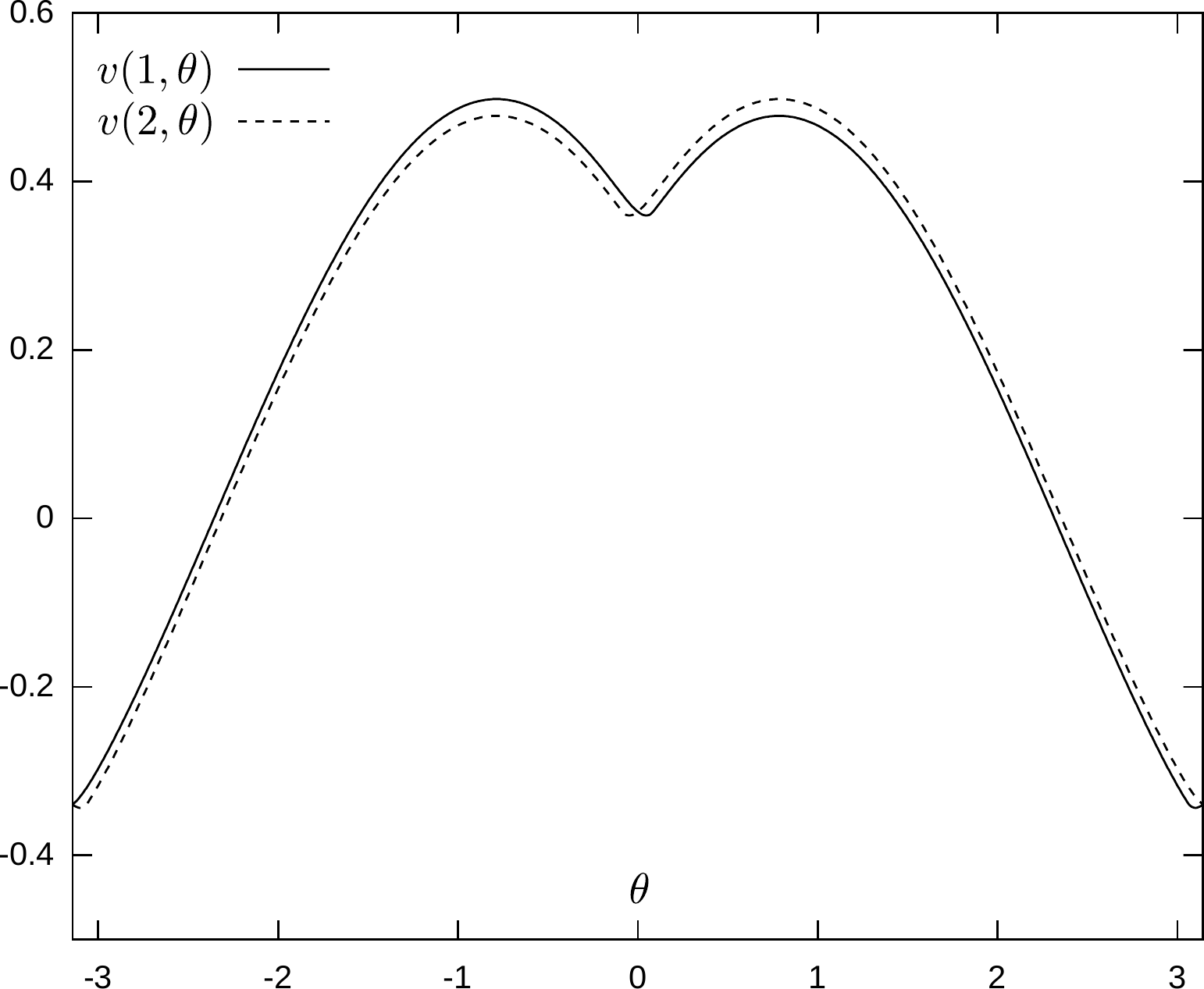}\\
(a)  
\end{tabular}
\vspace{0.3cm}
\begin{tabular}{cc}
\includegraphics[width=.45\textwidth]{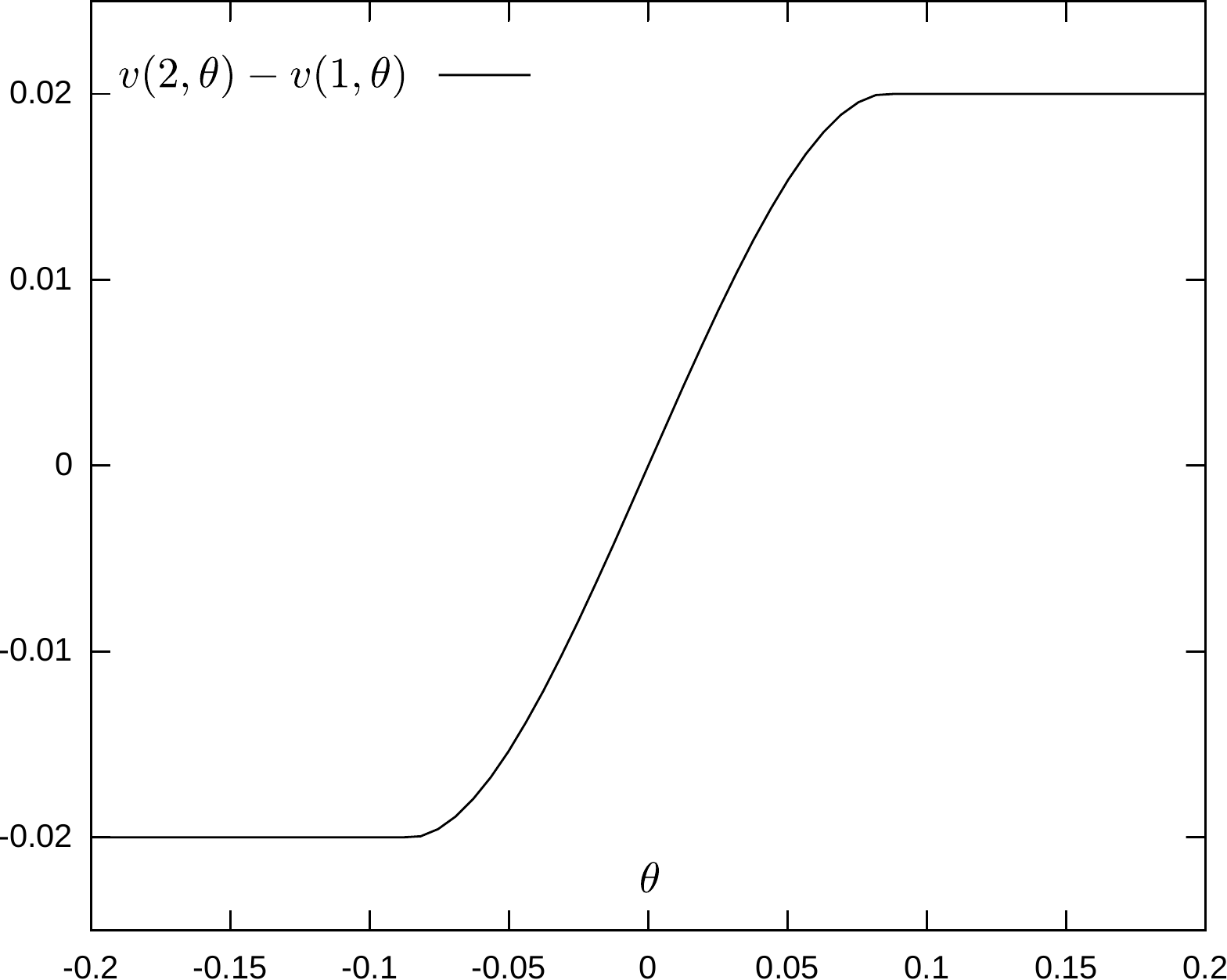}&
\includegraphics[width=.45\textwidth]{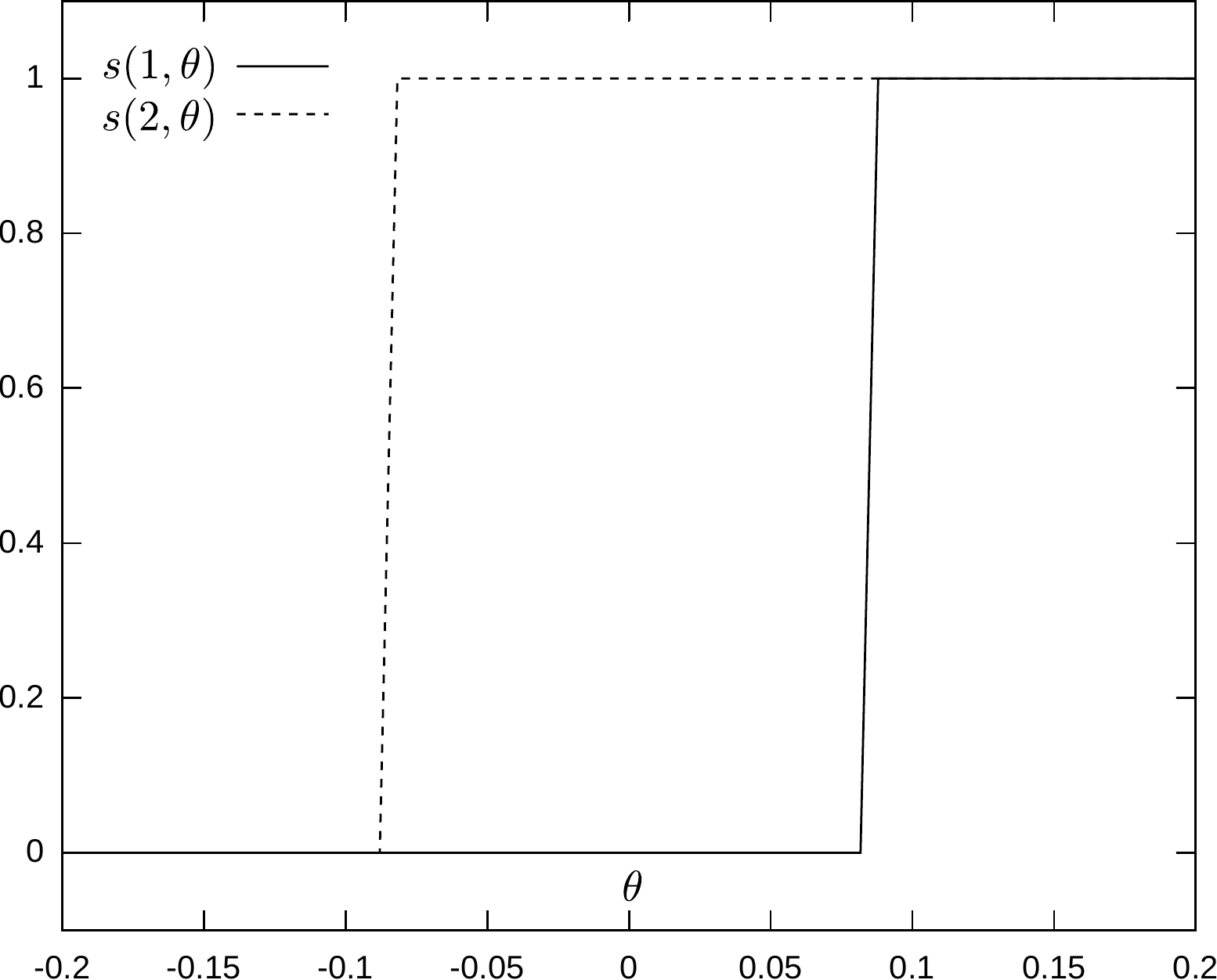}\\
(b)&(c) \\
\end{tabular}
\end{center}
\caption{One-dimensional problem. Value functions (a), zoom of their difference around the origin (b), and optimal switching maps (c).}\label{test0}
\end{figure}

\paragraph{\bf Example} We solve the one dimensional problem \eqref{singleplayerQVI-theta} for a single player, choosing the algorithm and parameters as described in the next section. Figure \ref{test0}a shows the value functions $v(1,\theta)$ and $v(2,\theta)$, corresponding to the two discrete states. We observe two crossing points, one at the origin, and one at the boundary of the periodic domain $[-\pi,\pi]$. In Figure \ref{test0}b we report, for $\theta\in[-0.2,0.2]$, a detail of the difference $v(2,\theta)-v(1,\theta)$, namely, the solution of the double obstacle problem \eqref{double-obstacle}. The computed contact point is about $\theta^*=0.085722$. Finally, in Figure \ref{test0}c, we show the optimal switching maps, observing a typical hysteresis loop around the origin, with optimal switching points $-\theta^*$ and $\theta^*$. 

\section{Numerical approximation}\label{numerics}
In this section we introduce a numerical scheme for solving the system of quasi-variational inequalities \eqref{mr-qvi}.
To this end, it is useful to rewrite \eqref{mr-qvi} in the following fixed point form:
\begin{equation}\label{mr-qvi-fixed-point}
v(q,r)=\max\left\{v(\hat q,r)-C^A,\min\left\{v(q,\hat r)+C^B, \frac{1}{\lambda}\left( f\cdot D v + \ell +\frac{\sigma^2}{2}\frac{\partial^2 v}{\partial \theta^2} \right) \right\}\right\}\,.
\end{equation}
Now, given $b_1,b_2,b_3>0$, we consider the computational box $[-b_1,b_1]\times[-b_2,b_2]\times[-b_3,b_3]$ in the reduced state space $\R^2\times[-\pi,\pi]$, 
and we introduce a uniform grid with nodes 
$$
\left(x_1^i,x_2^j,\theta^k\right)=(-b_1+i\Delta x_1,-b_2+j\Delta x_2,-b_3+k\Delta \theta),\qquad (i,j,k=0,\ldots,N),
$$
where $N$ is an integer and the space steps are given respectively by 
$\Delta x_1=2 b_1/N$, $\Delta x_2=2 b_2/N$ and $\Delta x_3=\Delta \theta=2 b_3/N$. \\
For a generic scalar or vector function $\chi(x_1,x_2,\theta,q,r)$, we denote by 
$\chi^{i,j,k}_{q,r}$ the corresponding approximation at the point $(x_1^i,x_2^j,\theta^k)$. Then, 
we discretize the stationary advection-diffusion equation appearing in \eqref{mr-qvi-fixed-point} 
using centred differences for the second derivative of $v$, and upwind differences (according to the sign of the components of $f$) for first derivatives. 
By straightforward algebraic manipulations, we obtain the following scheme:
\begin{equation}\label{fixpoint-scheme}
v^{i,j,k}_{q,r}=\mathcal{T}[v](i,j,k,q,r,\hat q,\hat r):=
\max\left\{v^{i,j,k}_{\hat q,r}-C^A,\min\left\{v^{i,j,k}_{q,\hat r}+C^B, \mathcal{S}[v^{i,j,k}_{q,r}] \right\}\right\}\,,
\end{equation}
where
\begin{equation}
\mathcal{S}[v^{i,j,k}_{q,r}]=\frac{1}{\Lambda}\left(\alpha_1 v^{\bar i,j,k}_{q,r}
                         +\alpha_2 v^{i,\bar j,k}_{q,r}
                         +\alpha_3 v^{i,j,\bar k}_{q,r}
                         +\frac12\alpha_4(v^{i,j,k-1}_{q,r}+v^{i,j,k+1}_{q,r})+l^{i,j,k}_{q,r}\right)
\end{equation}
with
$$
\alpha_1=\frac{|(f^{i,j,k}_{q,r})_1|}{\Delta x_1}\,,\quad\alpha_2=\frac{|(f^{i,j,k}_{q,r})_2|}{\Delta x_2}\,,\quad\alpha_3=\frac{|(f^{i,j,k}_{q,r})_3|}{\Delta \theta}
\,,\quad\alpha_4=\frac{\sigma^2}{\Delta \theta^2}\,,
$$
$$
\Lambda=\lambda +\alpha_1+\alpha_2+\alpha_3+\alpha_4
$$
and (the symbol $\sgn(\cdot)$ denotes the sign of its argument)
$$
\bar i=i+\sgn((f^{i,j,k}_{q,r})_1)\,,\quad\bar j=j+\sgn((f^{i,j,k}_{q,r})_2)\,,\quad\bar k=k+\sgn((f^{i,j,k}_{q,r})_3)\,.
$$
Now, we can compute the solution of \eqref{mr-qvi-fixed-point} using fixed point iterations, as described in Algorithm \ref{ALG1}.
\begin{algorithm}[ht]
\caption{Value Iteration Algorithm}\label{ALG1}
\begin{algorithmic}[1]
\STATE Assign an initial guess $(v^{i,j,k}_{q,r})^{(0)}$, for $i,j,k=0,\ldots,N$ and $q,r=1,2$. \\
Fix a tolerance $tol>0$ and set $n=0$\\
\REPEAT
\FOR {$i,j,k=1,\ldots,N-1$ and $q,r=1,2$}
\STATE Set $\hat q=3-q$ and $\hat r=3-r$
\STATE Compute $(v^{i,j,k}_{q,r})^{(n+1)}=\mathcal{T}[v^{(n)}](i,j,k,q,r,\hat q,\hat r)$
\ENDFOR
\STATE Set $n=n+1$
\UNTIL{$\displaystyle\max_{q,r}\max_{i,j,k}\left|(v^{i,j,k}_{q,r})^{(n)}-(v^{i,j,k}_{q,r})^{(n-1)}\right| < tol$}
\end{algorithmic}
\end{algorithm}

Note that, in this form, the scheme is consistent, monotone and $L^\infty$ stable (see the analysis in \cite{FerZid:2014}), and therefore convergent via the Barles--Souganidis theorem \cite{BarSou:1991}, in all cases in which a comparison principle holds.

We remark that the fixed point iterations are performed at the internal nodes of the grid, hence the choice of the boundary conditions for the initial guess is crucial. As discussed in the previous section, if the bounds $b_1$ and $b_2$ are large enough, the game at the boundary decouples in two optimal control problems, one for each player, both described by the same system of quasi-variational inequalities \eqref{singleplayerQVI-theta}, in the only state variable $\theta\in[-\pi,\pi]$. These one-dimensional problems can be solved again via fixed point iterations, using the following discretization of \eqref{singleplayerQVI-theta} for $P=A,B$ and $p=q,r\in\{1,2\}$:
$$
(v^P)^k_p=\max \left\{ (v^P)^k_{\hat p}-C^P,\mathcal S^1[(v^P)^k_p]\right\},
$$
with
$$
\mathcal S^1[(v^P)^k_p]=\left(\lambda+\frac{\sigma^2}{\Delta \theta^2}\right)^{-1}\left(\frac12 \frac{\sigma^2}{\Delta \theta^2} \Big((v^P)^{k-1}_p+(v^P)^{k+1}_p\Big)+\bar s^P\sin\left(\theta^k+\frac{\pi}{4}(-1)^p\right)\right)\,,
$$
and imposing periodic boundary conditions at $\theta=\pm\pi$. Once the solutions $v^A$ and $v^B$ are computed, we set the boundary values $v^{i,j,k}_{q,r}=(v^A)^k_q-(v^B)^k_r$ 
for $i=0$ or $i=N$ or $j=0$ or $j=N$ and $0\le k\le N$. 
Note that this relation can be used also in the internal nodes, to define a reasonable initial guess and save some iterations for the convergence of Algorithm \ref{ALG1}.

We finally remark that, in the special case $C^A=C^B$ and $\bar s^A=\bar s^B$, we can alternately solve the nonlinear equation \eqref{opt-comm} by a standard root-finding algorithm, and build the initial guess using the explicit expression \eqref{opt-val} for the difference $v^A-v^B$.

We proceed by discussing how to build optimal trajectories for the game. 
With the value function $v$ at hand, we have, by construction, the following inequalities for all $i,j,k=0,...,N$, all $q,r=1,2$ and $\hat q=3-q$, $\hat r = 3-r$
$$
v^{i,j,k}_{\hat q,r}-C^A\le v^{i,j,k}_{q,r}\le v^{i,j,k}_{q,\hat r}+C^B\,.
$$
Whenever an inequality is strict, the corresponding player keeps its discrete state, otherwise it can take an advantage on its opponent by switching to the other state and paying the corresponding cost. Then, we can easily define, for each player, an {\em optimal switching map}, depending on both the node $(x_1^i,x_2^j,\theta^k)$ and the state $(q,r)$: 
$$
{S^A}^{i,j,k}_{q,r}=\left\{
\begin{array}{ll}
 q &\mbox{if \,} v^{i,j,k}_{q,r}>v^{i,j,k}_{\hat q,r}-C^A \\
 \hat q&\mbox{if \,} v^{i,j,k}_{q,r}=v^{i,j,k}_{\hat q,r}-C^A
\end{array}
\right.
\qquad
{S^B}^{i,j,k}_{q,r}=\left\{
\begin{array}{ll}
 r &\mbox{if \,} v^{i,j,k}_{q,r}<v^{i,j,k}_{q,\hat r}+C^B \\
 \hat r&\mbox{if \,} v^{i,j,k}_{q,r}=v^{i,j,k}_{q,\hat r}+C^B
\end{array}
\right.
$$
Finally, we discretize the dynamics \eqref{game-dyn} by means of a simple forward Euler scheme with time step $\Delta t$:
$$
\left\{\begin{array}{l}
        X^A_{n+1}=X^A_n+f^A(X_n,\Theta_n,Q_n)\Delta t\\
        X^B_{n+1}=X^B_n+f^B(X_n,\Theta_n,R_n)\Delta t\\
        \Theta_{n+1}=\Theta_n+\sigma\sqrt{\Delta t} W_{n+1}\\
        Q_{n+1}={S^A}^{i_{n+1},j_{n+1},k_{n+1}}_{Q_n,R_n}\\
        R_{n+1}={S^B}^{i_{n+1},j_{n+1},k_{n+1}}_{Q_{n+1},R_n}
       \end{array}
\right.
\qquad
\left\{\begin{array}{l}
        X^A_0=x^A\\
        X^B_0=x^B\\
        \Theta_0=\theta\\
        Q_0=q\\
        R_0=r
       \end{array}
\right.
$$
where $\{W_n\}$ is a sequence of random numbers with a normal distribution of unit variance, and 
$$
\begin{array}{l}
i_{n+1}=\ceil{((X_{n+1})_1+b_1)/\Delta x_1}\\
j_{n+1}=\ceil{((X_{n+1})_2+b_2)/\Delta x_2}\\
k_{n+1}=\ceil{(\Theta_{n+1}+b_3)/\Delta \theta}\\
\end{array}
$$
define, by means of the upper integer part $\ceil{\cdot}$, a closest-neighbour projection on the grid of the updated state variables.

\section{Numerical examples}\label{examples}

Parameters for the simulations have been set according to the literature related to single-hull America's Cup vessels. In what follows, the length unit amounts to 1000 meters, and the time unit to 10 seconds. We choose the bounds $b_1=1$, $b_2=1$ and $b_3=\frac{\pi}{4}$, with $201$ nodes for each dimension of the grid (i.e., a total number of about $3.2\cdot 10^7$ nodes). 
Concerning the boat speeds, we choose $\bar s^A=\bar s^B=0.05$ and $\bar s^A_1=\bar s^B_1=300$. For the switching costs, we consider two different settings, a {\em symmetric} case with $C^A=C^B=0.02$, and 
an {\em asymmetric} case with $C^A=0.02$ and $C^B=0.04$. For the wind evolution, we consider a brownian motion with standard deviation 
$\sigma=0.03$. Finally, we set $\lambda=0.1$ for the discount factor in the cost functional, $tol=10^{-5}$ for the convergence tolerance in Algorithm \ref{ALG1}, and $\Delta t=0.2$ for the time step in the reconstruction of the optimal trajectories.

 As already remarked, uniqueness of solutions for the system \eqref{mr-qvi-fixed-point} is not ensured in the symmetric case $C^A=C^B$. Nevertheless, in the following tests, we always observe the convergence of the algorithm to a meaningful solution.  

\begin{figure}[!t]
\begin{center}
\begin{tabular}{cc}
\includegraphics[width=.475\textwidth]{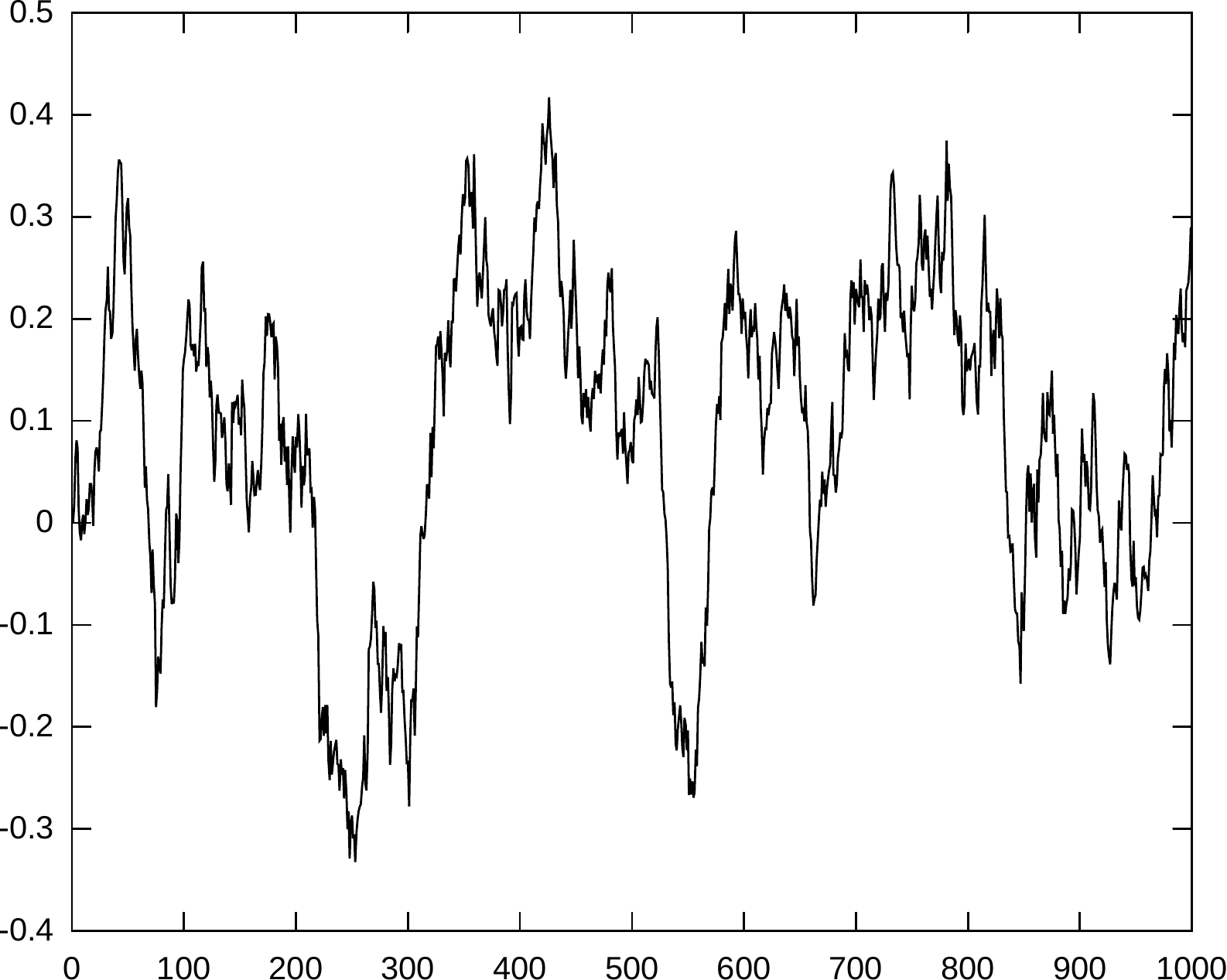}&
\includegraphics[width=.475\textwidth]{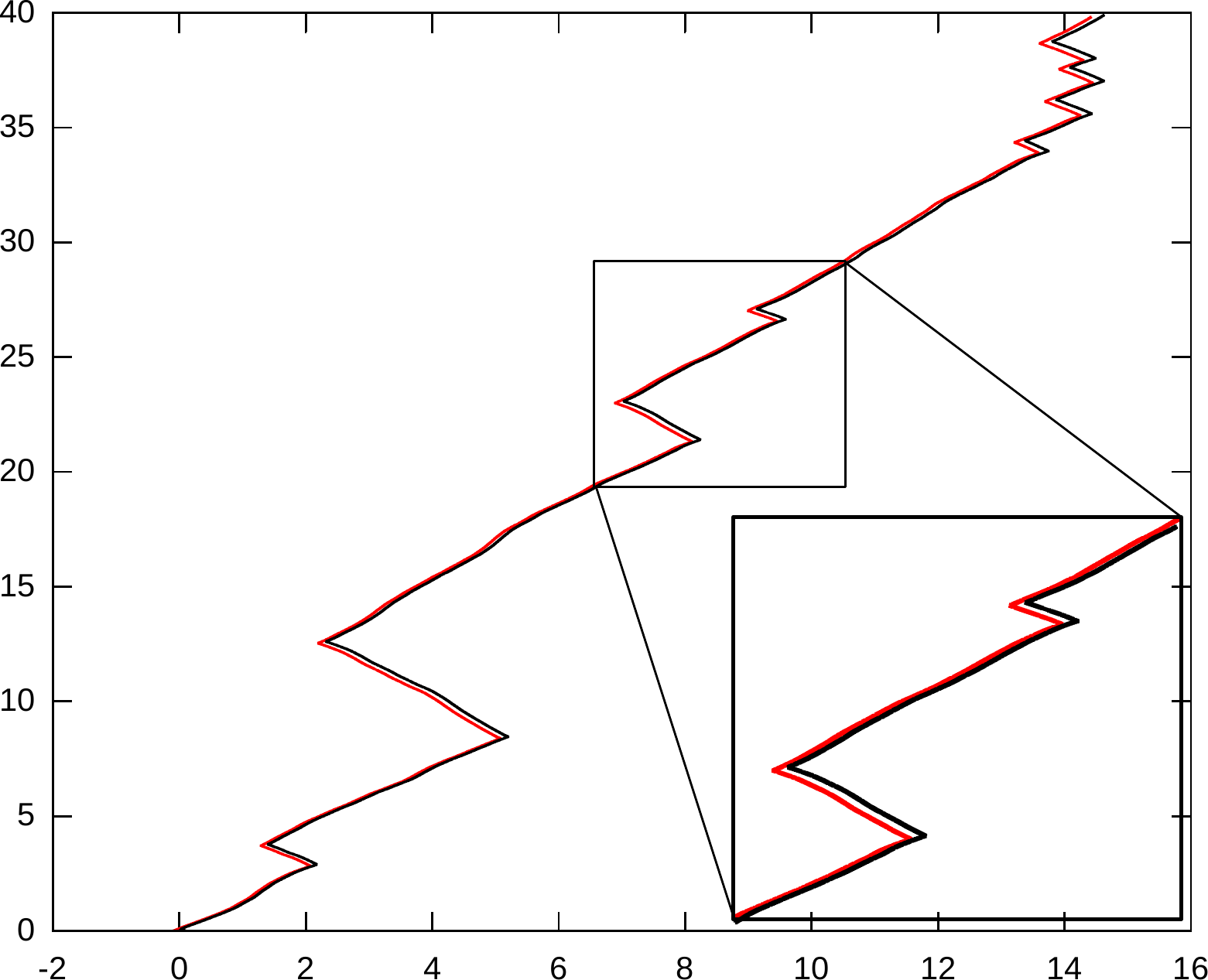}\\
wind direction & trajectories \\ \\
\includegraphics[width=.475\textwidth]{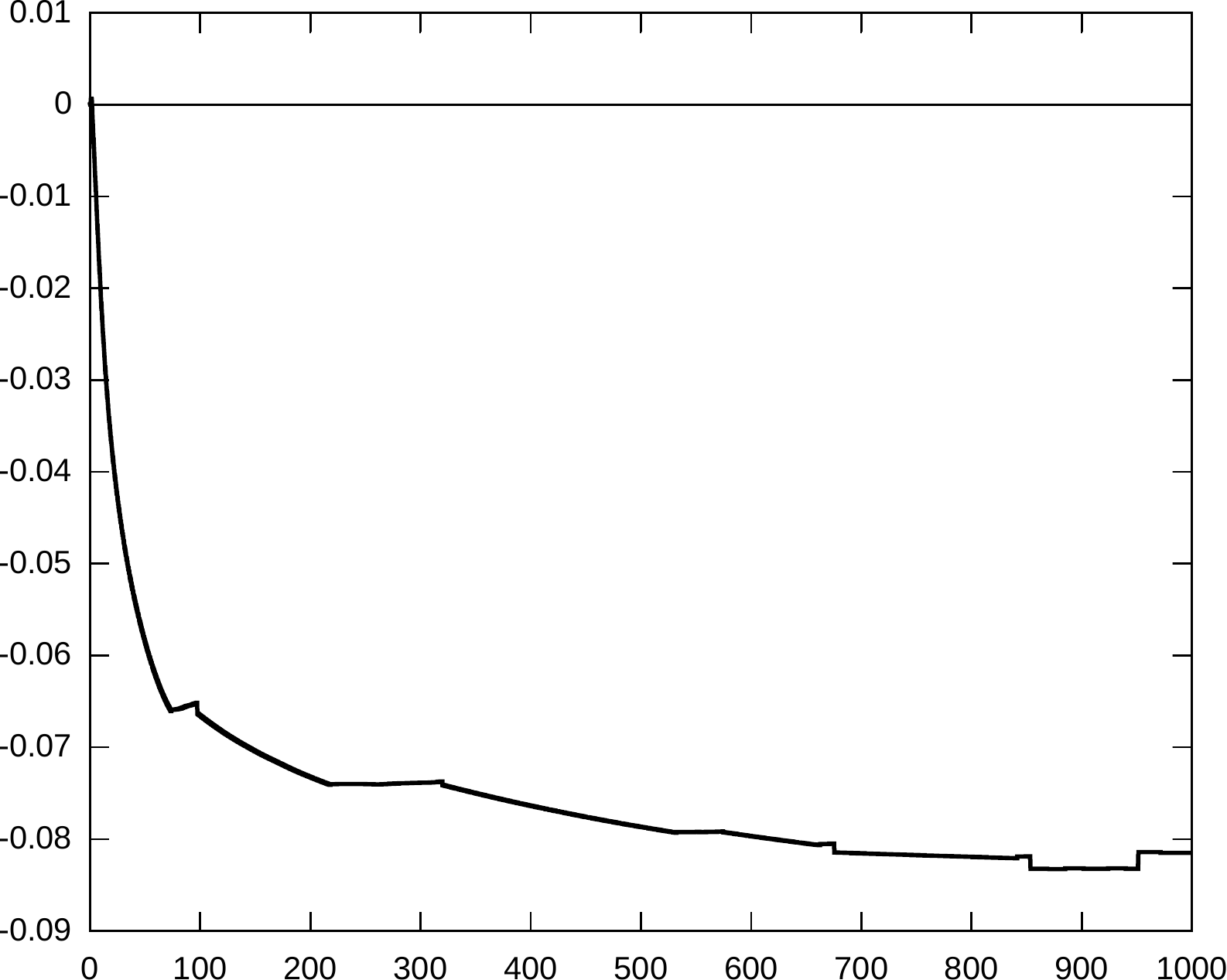}&
\includegraphics[width=.475\textwidth]{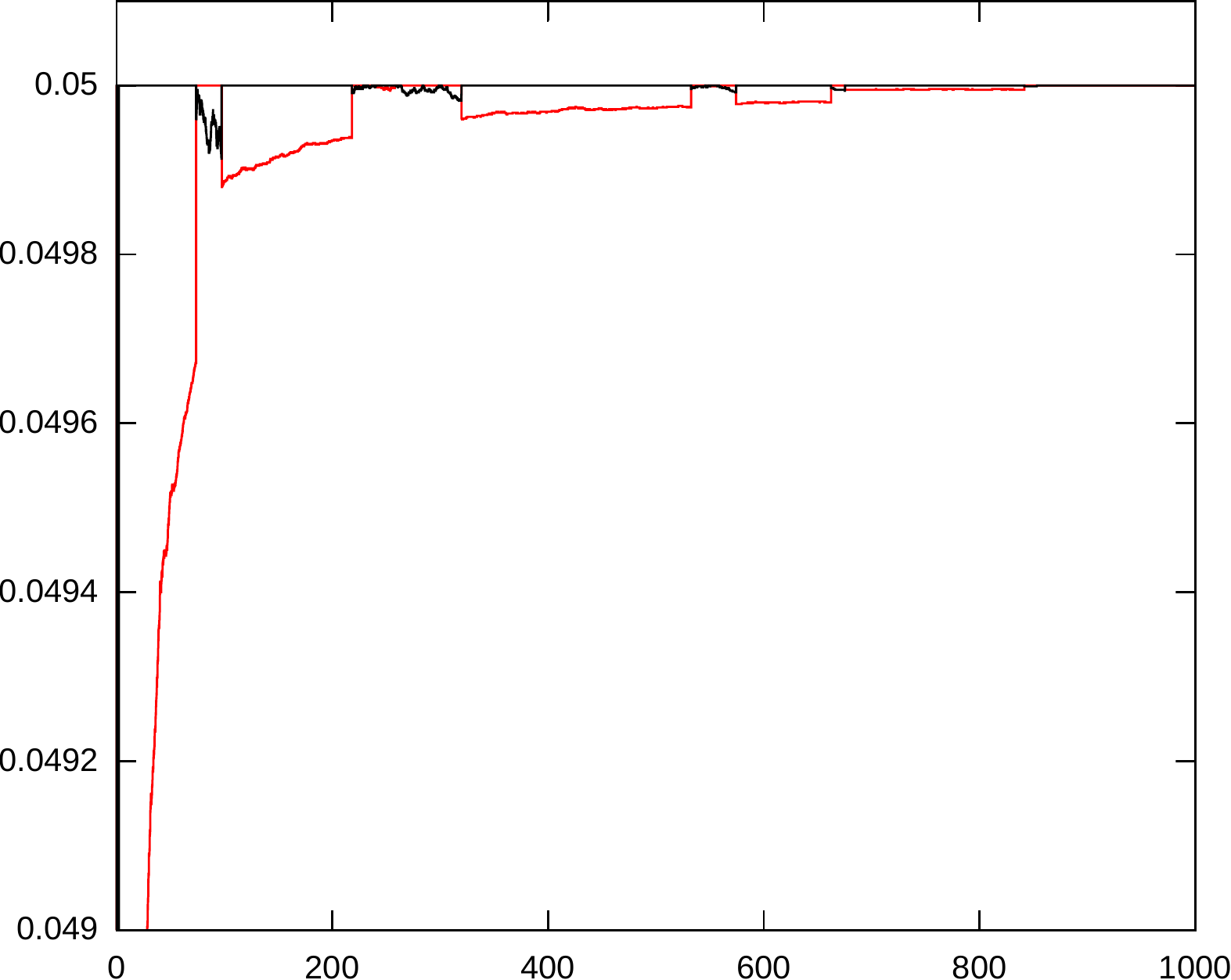} \\
 relative distance& speeds \\
\end{tabular}
\end{center}
\caption{Test 1a. Optimal strategy for both players in symmetric conditions, player B (black trajectory) wins.}\label{1a}
\end{figure}

\begin{figure}[!t]
\begin{center}
\begin{tabular}{cc}
\includegraphics[width=.475\textwidth]{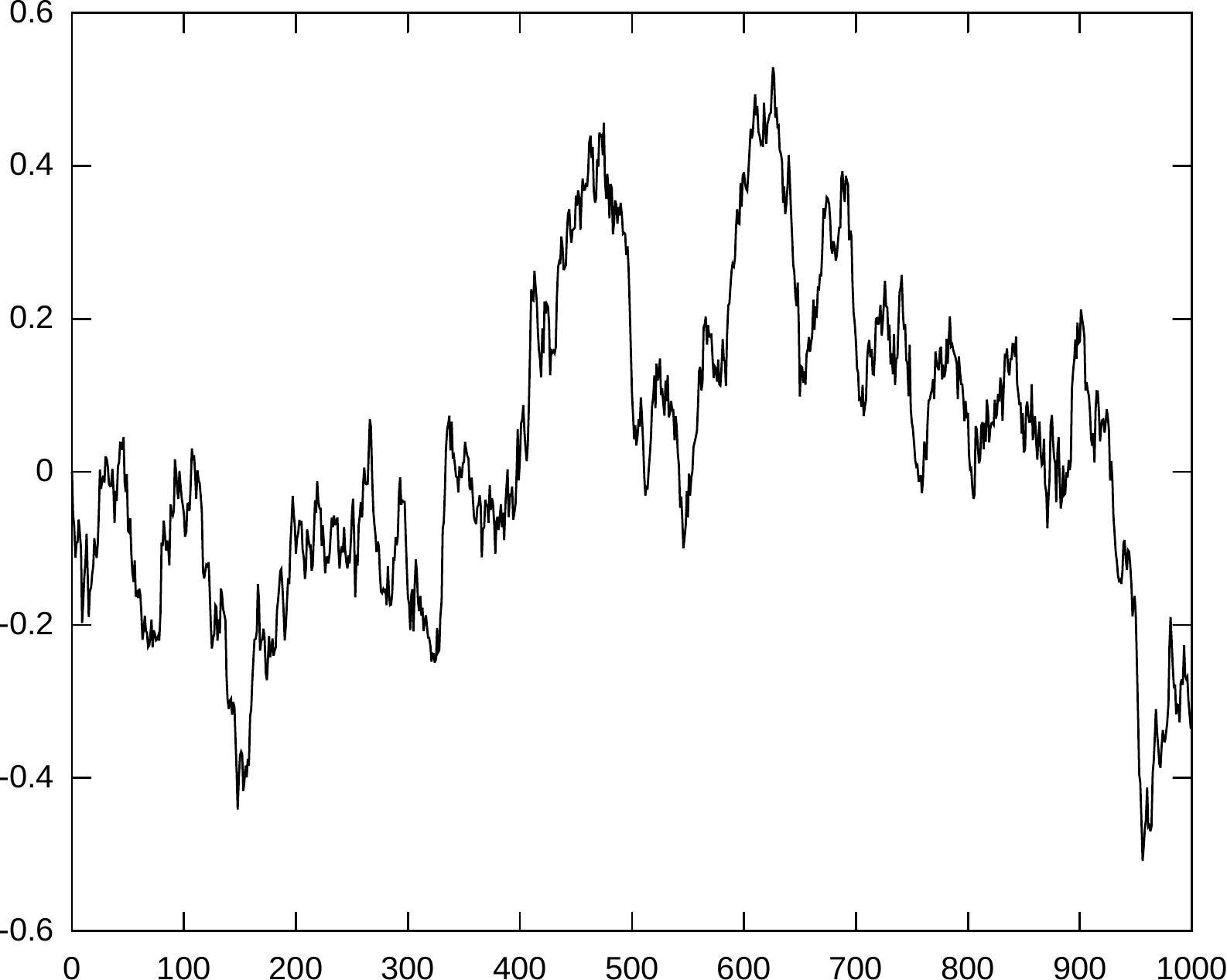}&
\includegraphics[width=.475\textwidth]{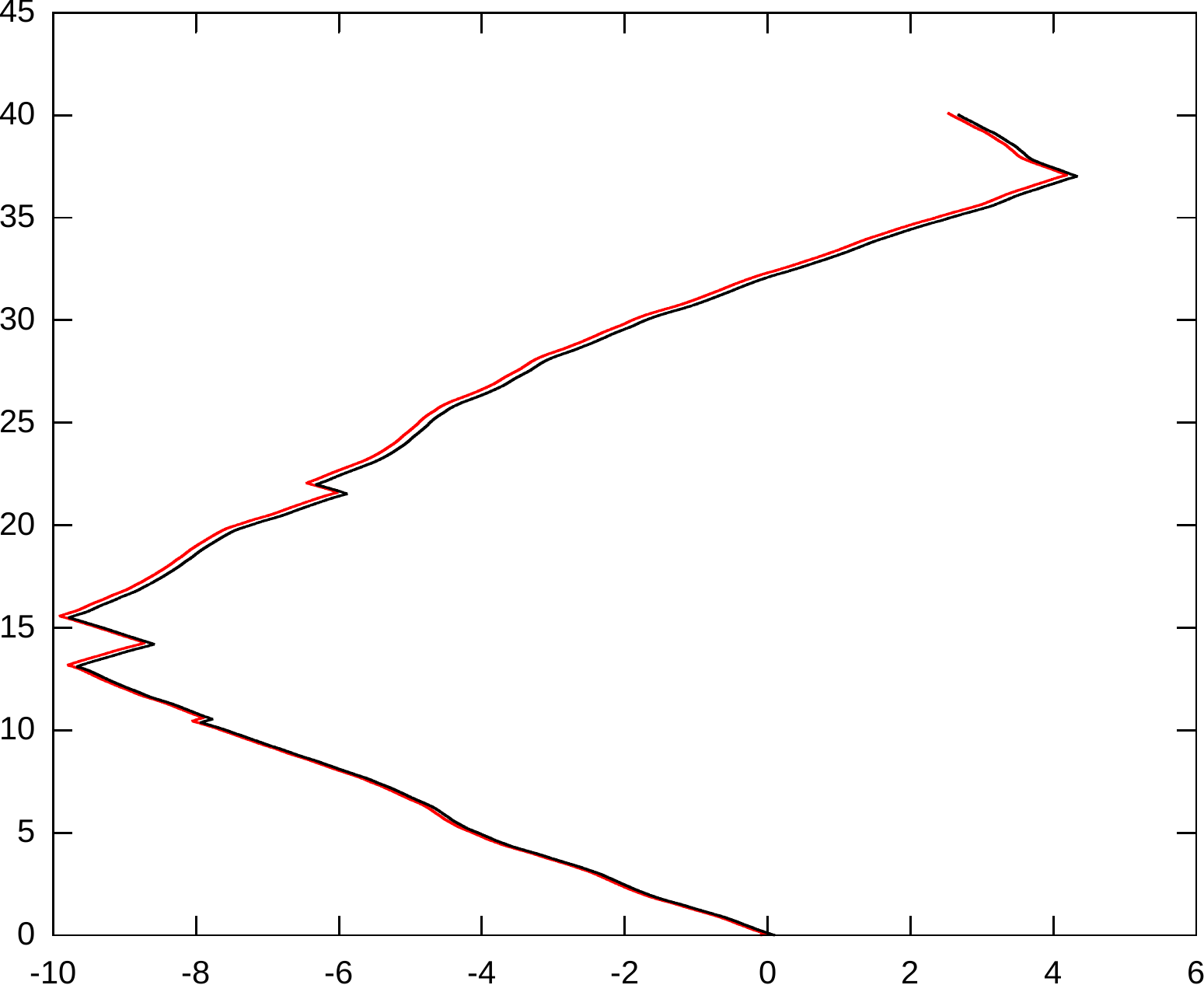}\\
wind direction& trajectories \\ \\
\includegraphics[width=.475\textwidth]{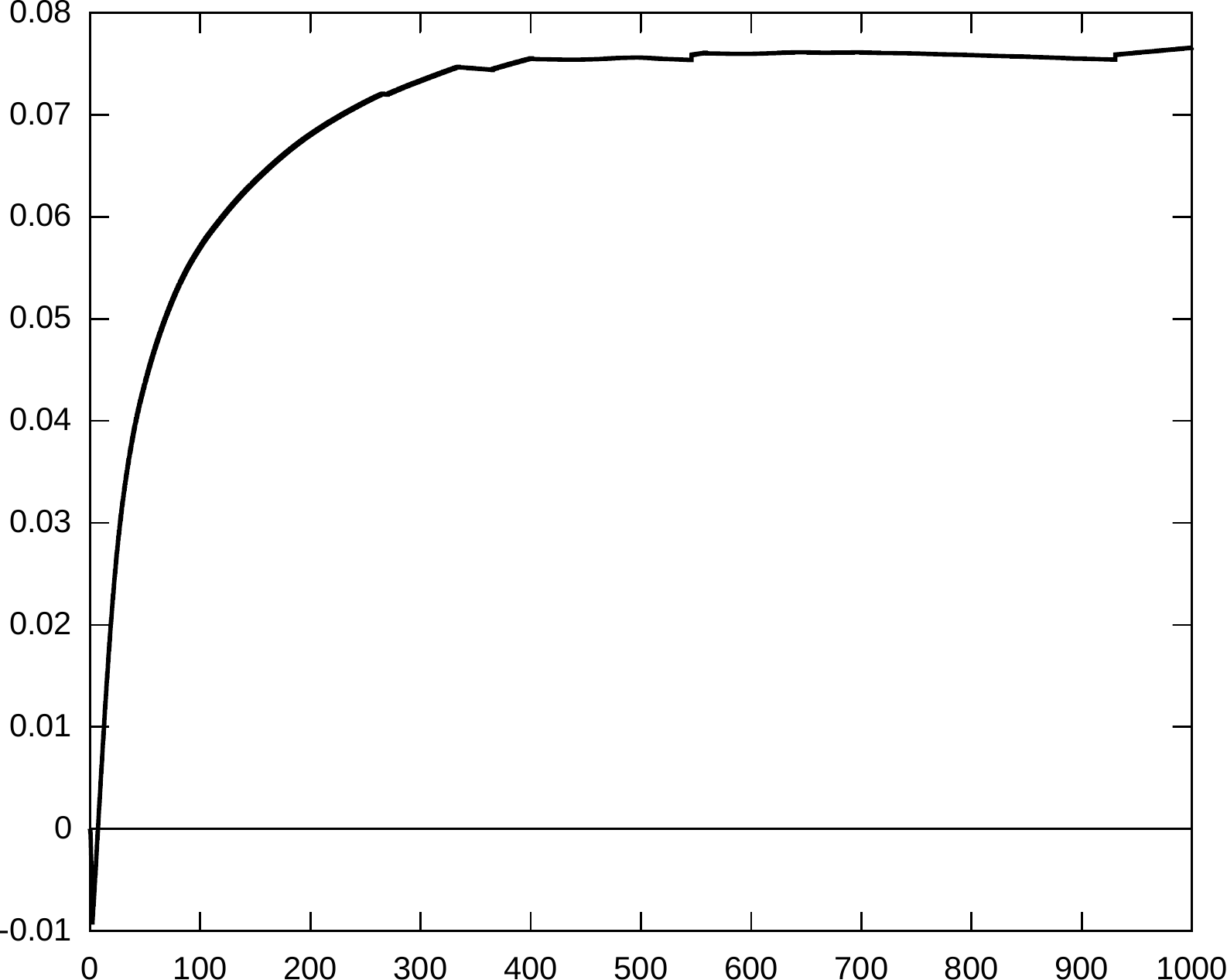}&
\includegraphics[width=.475\textwidth]{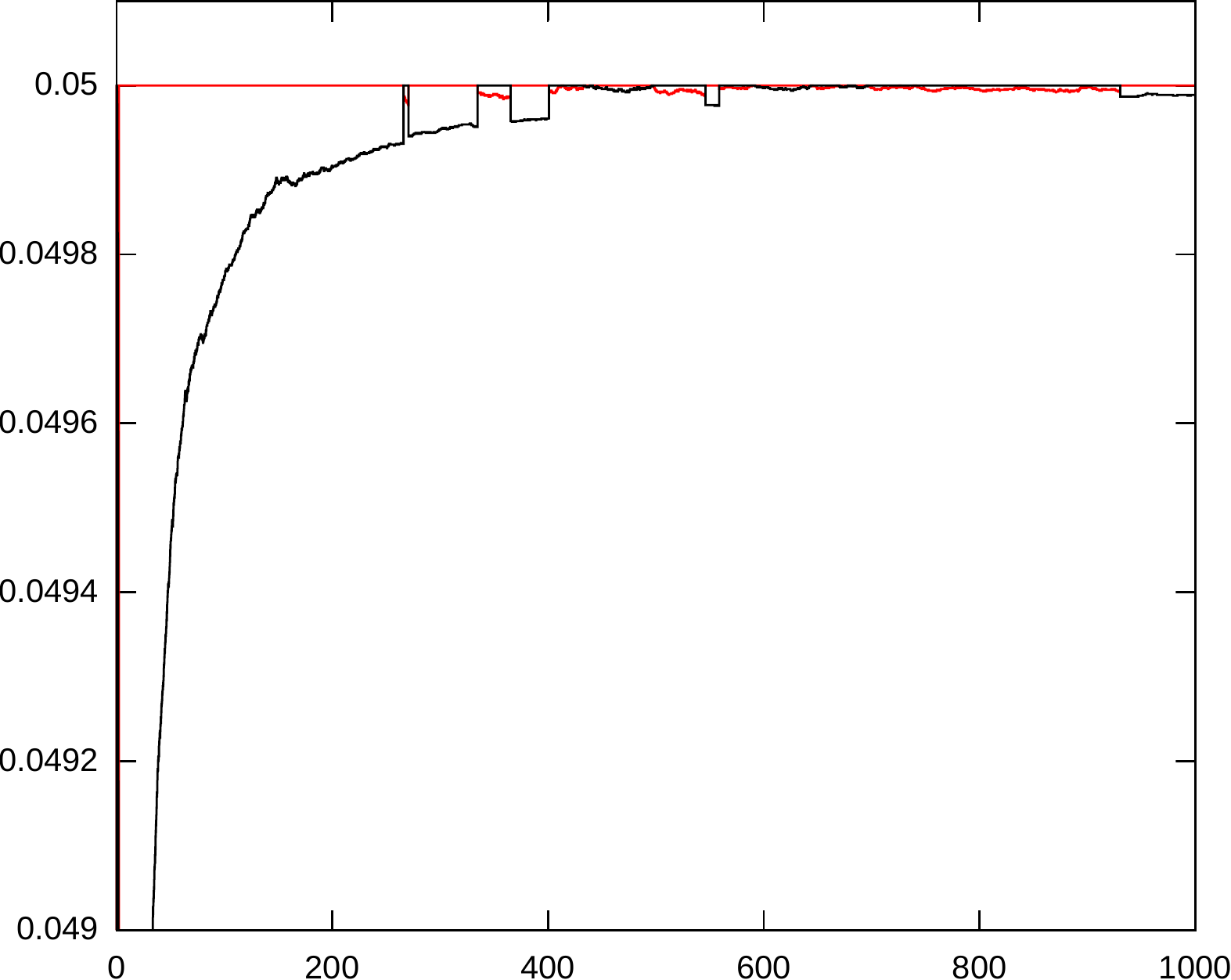} \\
relative distance& speeds \\
\end{tabular}
\end{center}
\caption{Test 1b. Optimal strategy for both players in symmetric conditions, player A (red trajectory) wins.}\label{1b}
\end{figure}

In the examples, we show some sample simulations obtained in typical scenarios. For each scenario, the value function and switching map have been computed in a first phase, while sample optimal (or suboptimal, as in the second example) trajectories are computed in the second phase, according to the procedure outlined in the previous section. For each simulation, four plots show respectively the wind evolution $\Theta(t)$, the resulting trajectories of the players in the $x_1-x_2$ plane, the relative position $x_2^A-x_2^B$ and the speeds of the two players, as functions of time. Trajectories and speeds are shown in red for player A, in black for player B.
\paragraph{\bf Test 1} We consider the symmetric case $C^A=C^B$, and the same initial $x_2$-coordinate, with the player A on the left side. Figures \ref{1a}--\ref{1b} show two sample trajectories.
The game is led, at least up to the final time $T=1000$ of the simulation, one time for each player. Both players tend to follow the optimal single-player strategy. However, the speed plots show that, once one of the players has gained a small advantage in the first part of the game, it tries to preserve the advantage by disturbing the other player as much as possible when in favourable position, and keeping away from the other if in unfavourable position. This results in two trajectories relatively close to one another, see also the detail of the trajectories in Fig. \ref{1a}.

\begin{figure}[!t]
\begin{center}
\begin{tabular}{cc}
\includegraphics[width=.475\textwidth]{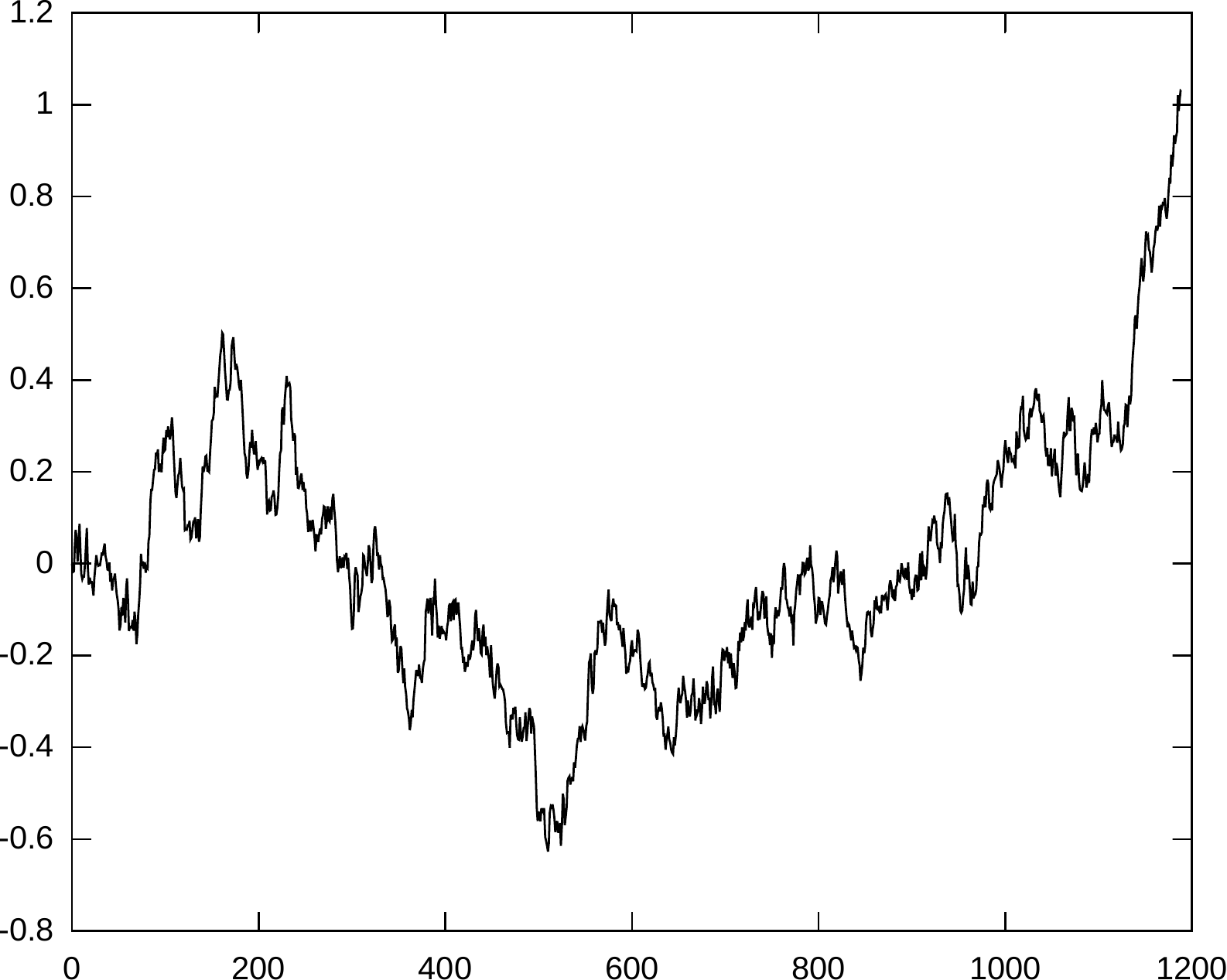}&
\includegraphics[width=.475\textwidth]{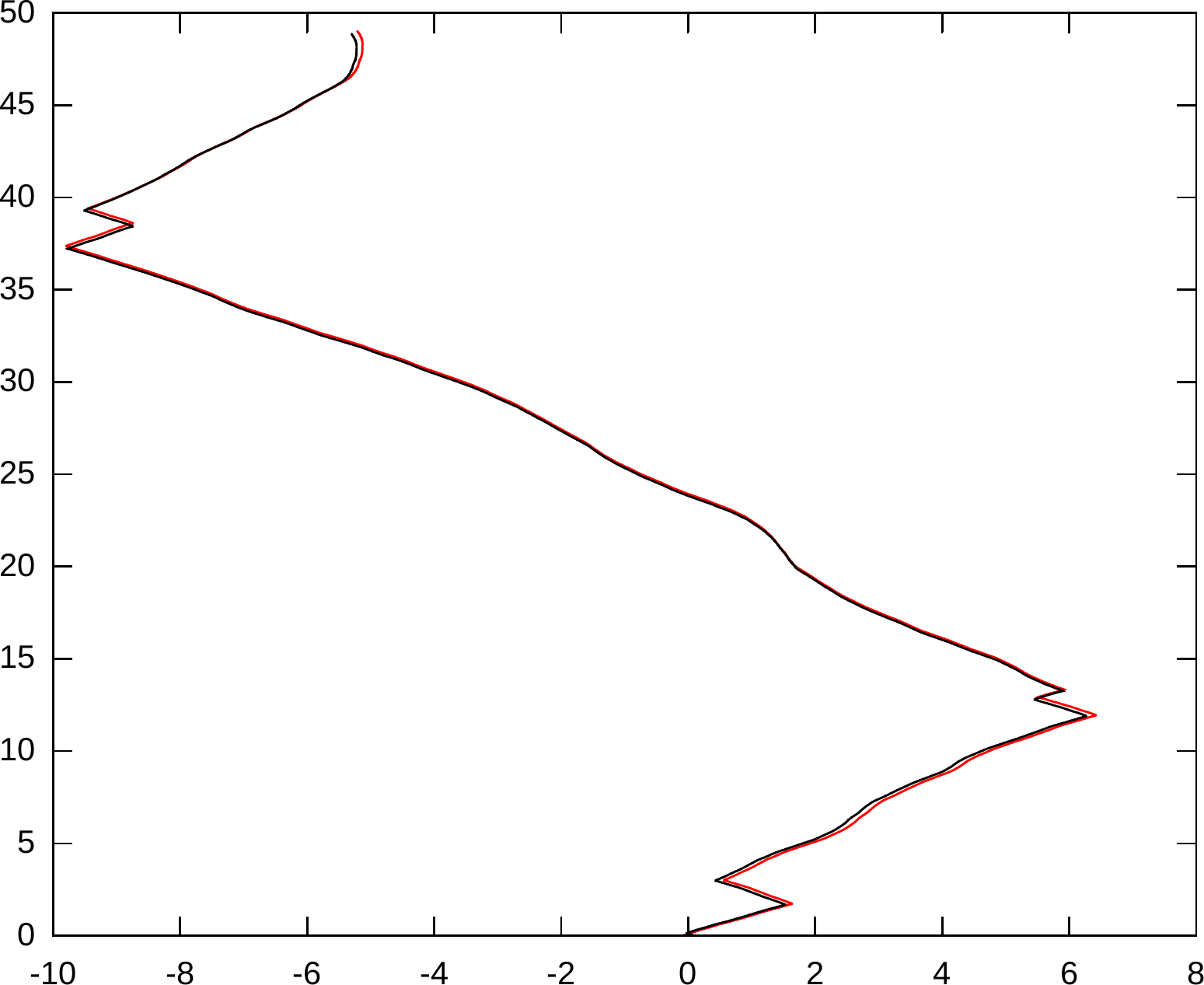}\\
wind direction& trajectories \\ \\
\includegraphics[width=.475\textwidth]{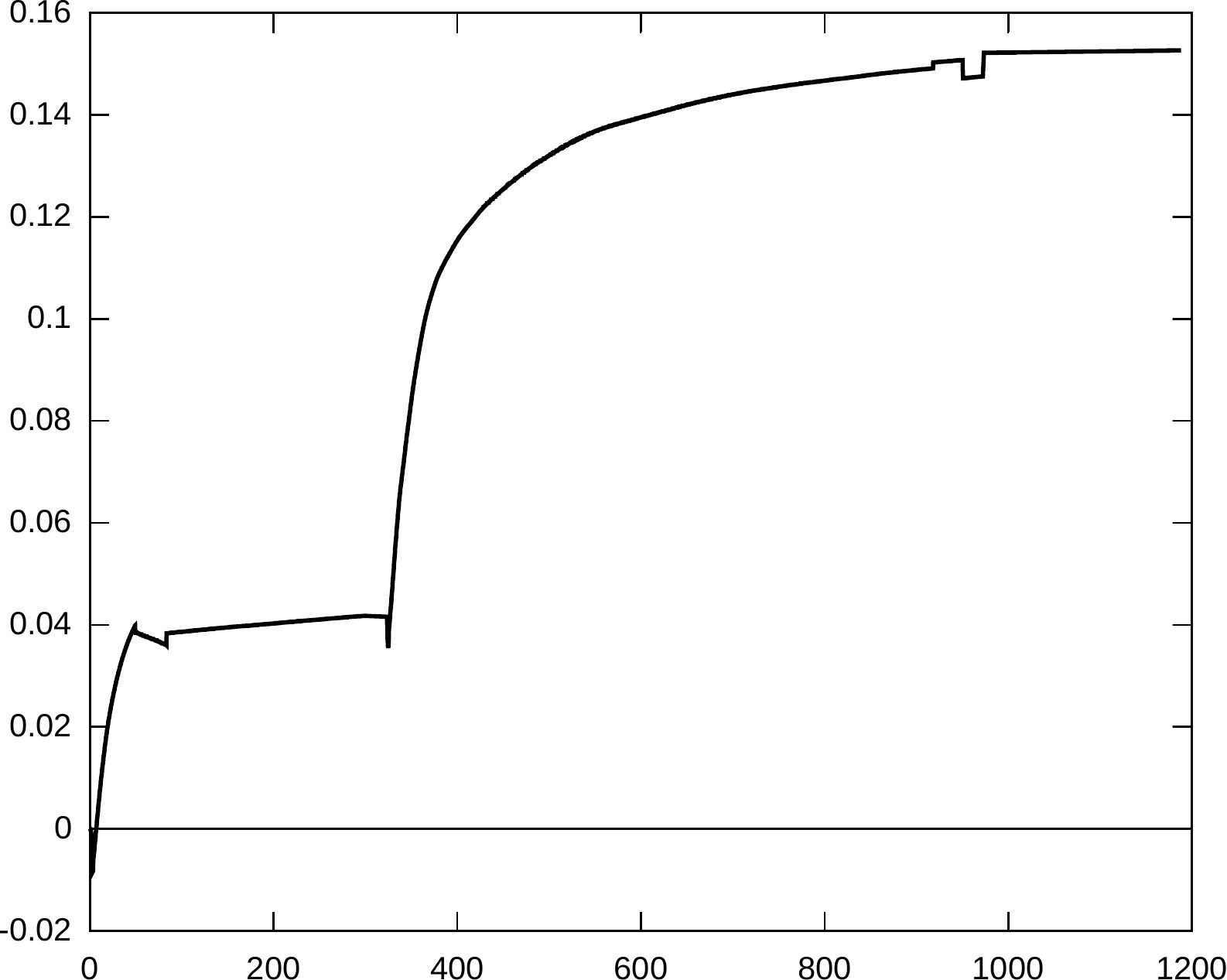}&
\includegraphics[width=.475\textwidth]{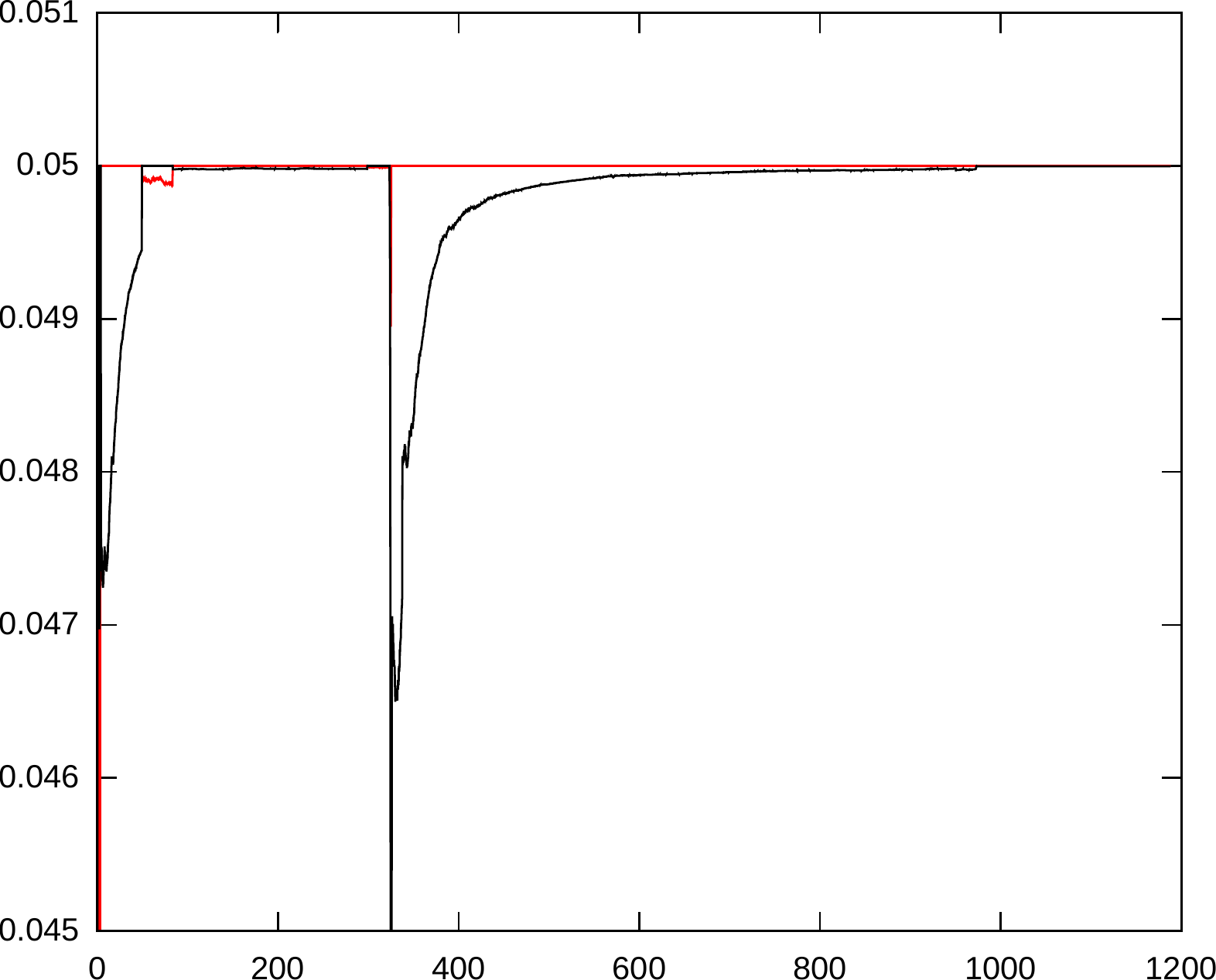} \\
relative distance& speeds \\
\end{tabular}
\end{center}
\caption{Test 2. Symmetric conditions, A (red) plays the optimal strategy for the two-player game, while B (black) plays the single-player optimal strategy.}\label{2}
\end{figure}

\paragraph{\bf Test 2} We still consider the symmetric case $C^A=C^B$, but here player A plays using the optimal strategy for the game, while player B plays using the single-player optimal strategy. Despite the small advantage gained by B in the first phase, A plays to disturb B (as it is apparent from the speed plot), and ends by leading the game, see Fig. \ref{2}.

\begin{figure}[!t]
\begin{center}
\begin{tabular}{cc}
\includegraphics[width=.475\textwidth]{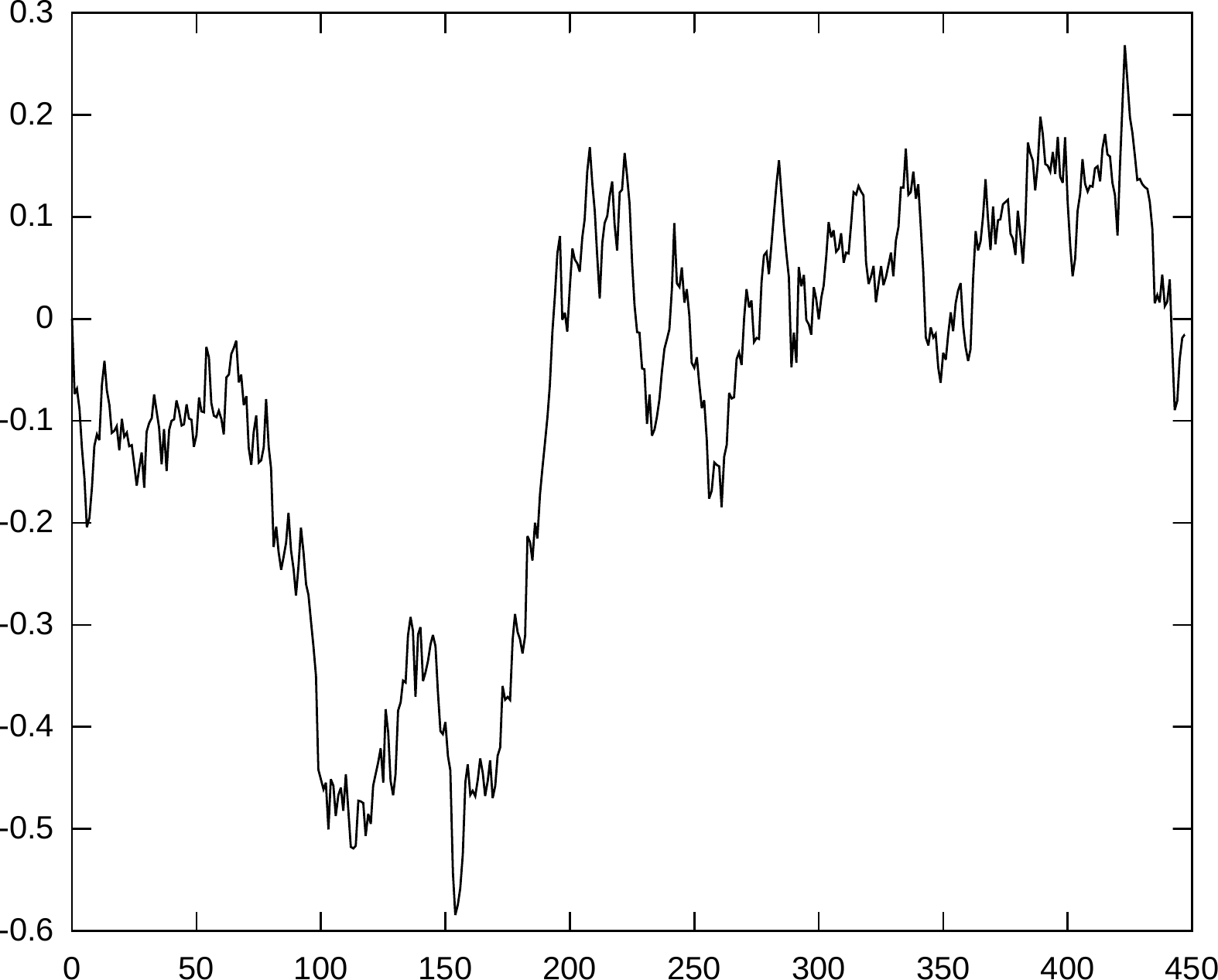}&
\includegraphics[width=.475\textwidth]{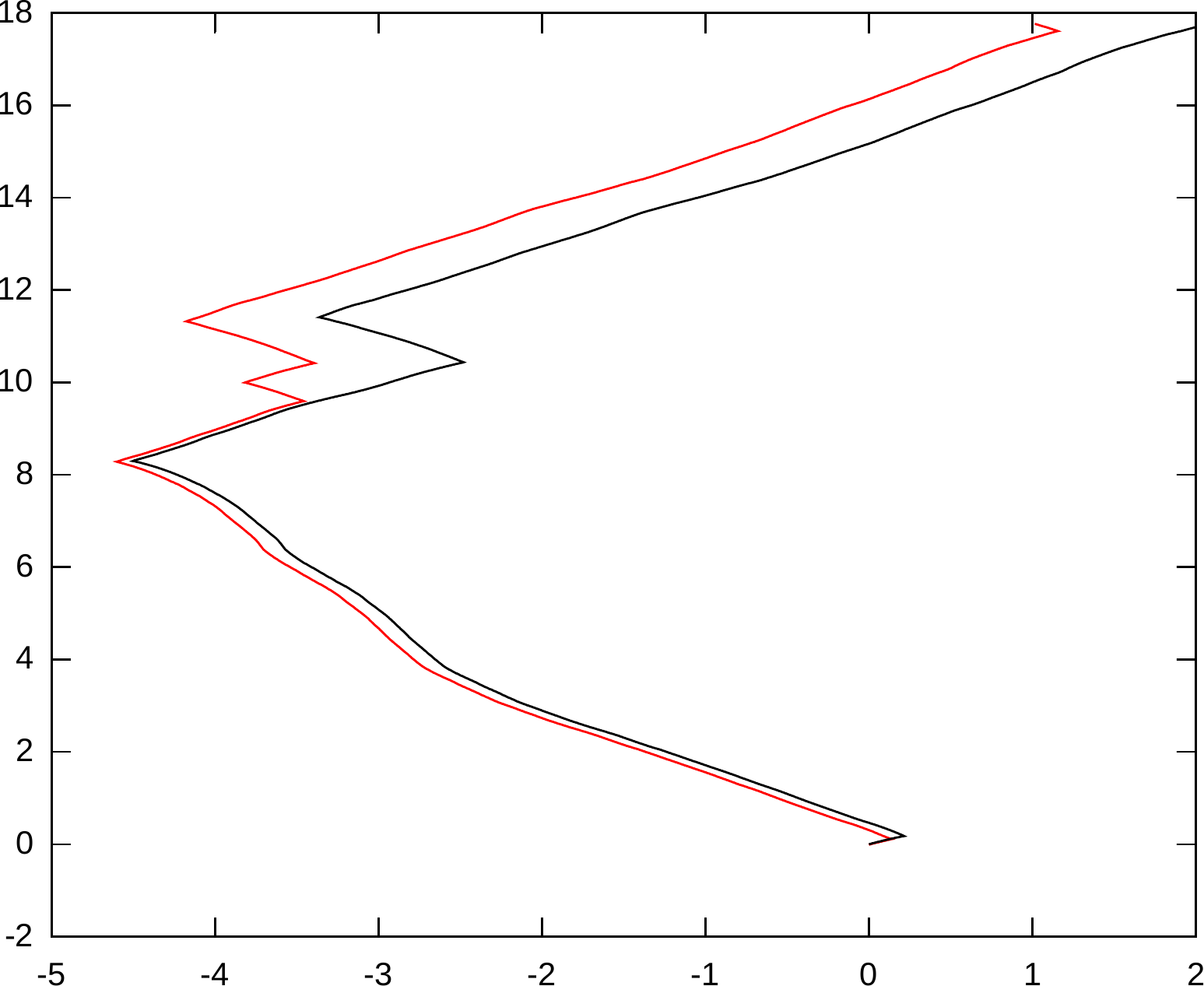}\\
wind direction& trajectories \\ \\
\includegraphics[width=.475\textwidth]{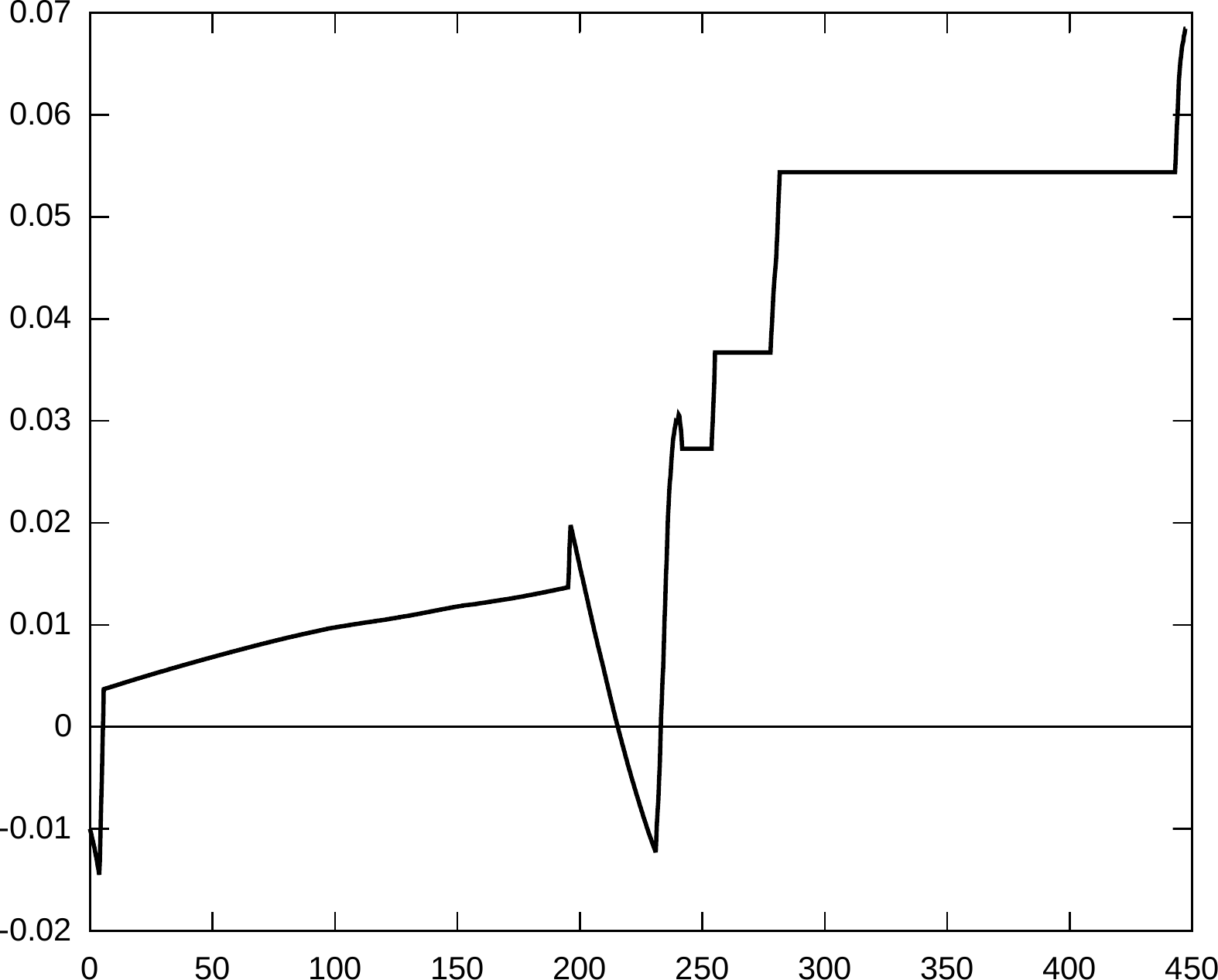}&
\includegraphics[width=.475\textwidth]{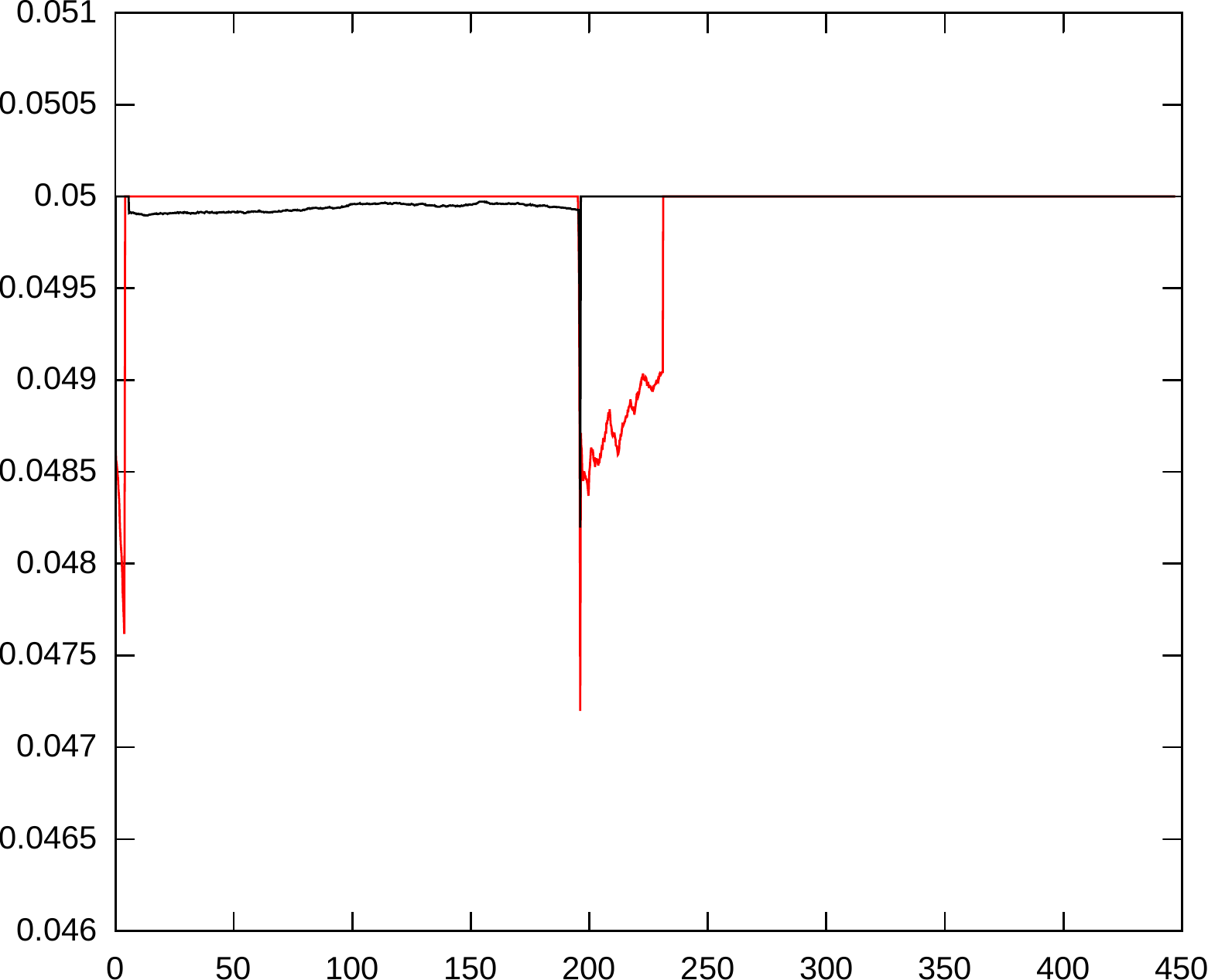} \\
relative distance& speeds \\
\end{tabular}
\end{center}
\caption{Test 3. Asymmetric conditions, player B (black) is ahead at the start, but pays a higher switching cost.}\label{3}
\end{figure}

\paragraph{\bf Test 3} We finally consider the asymmetric case $C^A\neq C^B$. Player B starts in a favourable position, but player A ends by leading the game. Here, rather than from the coupling between the players, A seems to take advantage of its better ability to exploit wind variations, see Fig. \ref{3}.

\section*{References}
\bibliographystyle{elsarticle-num}
\bibliography{references}             

\end{document}